\documentclass[11pt,reqno,a4paper]{amsart}

\usepackage{amssymb,mathtools,enumitem,microtype}
\usepackage[colorlinks,linkcolor=magenta,anchorcolor=blue,citecolor=green]{hyperref}

\hypersetup{pdftitle={Almost all vertex-transitive graphs and digraphs are Cayley},pdfauthor={Ning Su and Binzhou Xia}}

\newtheorem{theorem}{Theorem}[section]
\newtheorem{proposition}[theorem]{Proposition}
\newtheorem{lemma}[theorem]{Lemma}

\theoremstyle{definition}
\newtheorem{definition}[theorem]{Definition}
\newtheorem{remark}[theorem]{Remark}

\newcommand{\AGL}{\mathrm{AGL}}
\newcommand{\Alt}{\mathrm{Alt}}
\newcommand{\Aut}{\mathrm{Aut}}
\newcommand{\Cay}{\mathrm{Cay}}
\newcommand{\core}{\mathrm{Core}}
\newcommand{\End}{\mathrm{End}}
\newcommand{\GL}{\mathrm{GL}}
\newcommand{\Ind}{\mathrm{Ind}}
\newcommand{\Inn}{\mathrm{Inn}}
\newcommand{\Ker}{\mathrm{Ker}}
\newcommand{\lpp}{\mathrm{lpp}}
\newcommand{\Out}{\mathrm{Out}}
\newcommand{\Res}{\mathrm{Res}}
\newcommand{\Sym}{\mathrm{Sym}}

\title[Vertex-transitive graphs and digraphs]{Almost all vertex-transitive graphs and digraphs are Cayley}
\author{Ning Su}
\author{Binzhou Xia}
\date{}

\begin{document}

\begin{abstract}
We prove that almost all vertex-transitive digraphs are Cayley digraphs, and that almost all vertex-transitive graphs are Cayley graphs, as the number of vertices tends to infinity through all positive integers. The proportions of Cayley objects are at least $1-2^{-n/4+o(n)}$ for digraphs and $1-2^{-n/24+o(n)}$ for graphs. This proves the McKay--Praeger conjecture and its digraph analogue.
\end{abstract}

\maketitle

\section{Introduction}\label{sec:introduction}

In this paper, a digraph $\Gamma$ consists of a finite vertex set $V(\Gamma)$ and an arc set $A(\Gamma)$ that is an irreflexive binary relation on $V(\Gamma)$. A graph is a digraph whose arc set is symmetric.

A digraph $\Gamma$ is \emph{vertex-transitive} if $\Aut(\Gamma)$ is transitive on its vertex set.
For a finite group $R$ and a subset $S\subseteq R\setminus\{1\}$, the \emph{Cayley digraph} $\Cay(R,S)$ on $R$ has vertex set $R$ and arc set $\{(x,sx)\mid x\in R,\ s\in S\}$. It is a graph if and only if $S=S^{-1}$.
A digraph is called \emph{Cayley} if it is isomorphic to $\Cay(R,S)$ for some group $R$ and subset $S\subseteq R\setminus\{1\}$. It is well known that a digraph is Cayley if and only if its automorphism group contains a regular subgroup.
Every Cayley digraph is vertex-transitive, but the converse is not true.
The construction and classification of non-Cayley vertex-transitive graphs, including the determination of their possible orders, have been studied extensively; see, for example,~\cite{DS,IP,MP,MP2}.

For small orders, McKay and Praeger~\cite[p.~54]{MP} tabulated the numbers of vertex-transitive graphs and of those that are not Cayley.
Their counts are up to isomorphism at each fixed order. Observing that the great majority of these graphs are Cayley graphs, they expressed the expectation that this trend would continue at larger orders. This expectation has since been formulated as the McKay--Praeger conjecture: the proportion of Cayley graphs among vertex-transitive graphs tends to one as the order tends to infinity.
Our main result proves this conjecture with an exponentially small proportion of exceptions and establishes the corresponding extension to digraphs.

For a positive integer $n$, let $\mathcal V(n)$ and $\mathcal C(n)$ be the sets of isomorphism classes of vertex-transitive digraphs and Cayley digraphs of order $n$, respectively.
Let $\mathcal U(n)\subseteq\mathcal V(n)$ consist of the vertex-transitive graphs of order $n$.
Thus $\mathcal C(n)\cap\mathcal U(n)$ is the set of isomorphism classes of Cayley graphs of order $n$.

\begin{theorem}\label{thm:main}
As the positive integer $n\to\infty$,
\[
\frac{|\mathcal C(n)|}{|\mathcal V(n)|}\geq1-2^{-\frac{n}{4}+o(n)}\quad\text{and}\quad
\frac{|\mathcal C(n)\cap\mathcal U(n)|}{|\mathcal U(n)|}\geq1-2^{-\frac{n}{24}+o(n)}.
\]
\end{theorem}

Our proof uses the classification of finite simple groups.
We first count orbital partitions of minimally transitive groups, and then count invariant graphs and digraphs on representatives of these partitions.
The proof of Theorem~\ref{thm:main} is given in Section~\ref{sec:counting}.
% In Section~\ref{sec:counting} we also give the connected versions of Theorem~\ref{thm:main} in Corollary~\ref{cor:connected}.
% Moreover, Corollaries~\ref{cor:cumulative} and \ref{cor:connected-cumulative} give the cumulative versions, counting objects of order at most $n$, without and with connectedness, respectively.
% Corollary~\ref{cor:labelled} gives all four versions for labelled objects, with the same lower bounds on the proportions of Cayley objects.

\section{Preliminaries}\label{sec:prelim}

In this paper, all groups are finite, and all logarithms are to base $2$.
For a positive integer $n$ and a prime $p$, let $n_p$ denote the largest power of $p$ dividing $n$.
Group actions are on the right, and mappings are composed from left to right.
We write $x^g$ for the image of $x$ under $g$, and $\Inn(g)$ for the automorphism $x\mapsto g^{-1}xg$.

Let $G$ be a group. Denote by $\mathbf{Z}(G)$ the centre of $G$, by $d(G)$ the smallest size of a generating set of $G$, and by $i(G)$ the number of elements $x\in G$ such that $x^2=1$.
We write $G/H=\{Hg\mid g\in G\}$ for the set of right cosets of a subgroup $H$ of $G$.
The action of $G$ on $G/H$ by right multiplication has kernel $\core_G(H)=\bigcap_{g\in G}g^{-1}Hg$.

For a permutation group $G$ on a set $\Omega$, an \emph{orbital} is an orbit of $G$ on $\Omega\times\Omega$. Let $r(G)$ denote the number of $G$-orbitals.
The partition of $\Omega\times\Omega$ into $G$-orbitals is the \emph{orbital partition} of $G$.
For convenience, we introduce the following definition.

\begin{definition}\label{def:orbital-isomorphism}
Let $G\leq\Sym(\Omega)$ and $H\leq\Sym(\Delta)$.
An \emph{isomorphism} from the orbital partition of $G$ to that of $H$ is a bijection $f\colon\Omega\to\Delta$ such that $(u,v)\mapsto(u^f,v^f)$ maps each $G$-orbital onto an $H$-orbital.
The orbital partitions are called \emph{isomorphic} if such a bijection exists.
\end{definition}

Let $V$ be a finite-dimensional complex $G$-module. A $G$-endomorphism of $V$ is a linear mapping $\varphi\colon V\to V$ such that $\varphi(x^g)=\varphi(x)^g$ for all $x\in V$ and $g\in G$. We denote the algebra of these endomorphisms by $\End_G(V)$. The kernel of the action of $G$ on $V$ is the subgroup consisting of the elements that fix every vector of $V$. Isomorphic $G$-modules have the same kernel, since a $G$-module isomorphism carries the equality $x^g=x$ on one module to the same equality on the other.
For each $D\trianglelefteq G$, let $V_D$ be the sum of all irreducible $G$-submodules of $V$ with kernel $D$. In particular, $V_G$ consists of all vectors fixed by every element of $G$. Clearly,
\[
V=\bigoplus_{D\trianglelefteq G}V_D.
\]
Moreover, every $\varphi\in\End_G(V)$ preserves each $V_D$: on an irreducible submodule $U$, the restriction of $\varphi$ is either zero or an isomorphism onto its image; in the latter case, $U$ and its image have the same kernel.

For a finite $G$-set $\Omega$, the permutation module $\mathbb C[\Omega]$ has basis $\{e_\alpha\mid\alpha\in\Omega\}$, with $e_\alpha^g=e_{\alpha^g}$. Its subspace fixed pointwise by $G$ consists of the vectors whose coefficients are constant on each $G$-orbit on $\Omega$. The following lemma relates its endomorphisms to the orbits on ordered pairs.

\begin{lemma}\label{lem:orbital-endomorphisms}
Let $G\leq\Sym(\Omega)$ with $\Omega$ finite. The indicator matrices of the $G$-orbitals form a basis of $\End_G(\mathbb C[\Omega])$. Moreover, a permutation $s\in\Sym(\Omega)$ stabilizes every $G$-orbital if and only if its linear action on $\mathbb C[\Omega]$ commutes with every member of $\End_G(\mathbb C[\Omega])$.
\end{lemma}

\begin{proof}
Let $\varphi$ be a linear endomorphism of $\mathbb C[\Omega]$, and write $\varphi(e_\alpha)=\sum_{\beta\in\Omega}a_{\alpha,\beta}e_\beta$ for each $\alpha\in\Omega$. For $g\in G$, the equality $\varphi(e_\alpha^g)=\varphi(e_\alpha)^g$ for every $\alpha\in\Omega$ is equivalent to
\[
a_{\alpha^g,\beta^g}=a_{\alpha,\beta}\quad\text{for all }\alpha,\beta\in\Omega.
\]
Thus $\varphi\in\End_G(\mathbb C[\Omega])$ if and only if its matrix entries are constant on each $G$-orbital. The indicator matrices of these orbitals therefore form a basis. The same calculation with $s$ in place of $g$ shows that $s$ commutes with the indicator matrix of an orbital if and only if it stabilizes that orbital, proving the second assertion.
\end{proof}

A transitive permutation group is \emph{minimally transitive} if it has no proper transitive subgroup.
Parts~\ref{item:rank} and~\ref{item:generators} of the following lemma are~\cite[Lemma~3.1]{MS} and~\cite[Lemma~3.1]{NVL}, respectively, and Part~\ref{item:normalizer} follows from~\cite[Theorem~1.7]{GMP}.

\begin{lemma}\label{lem:rank-generators}
Let $G$ be a transitive permutation group of degree $n\geq1$.
The following statements hold:
\begin{enumerate}[label=\textup{(\alph*)},leftmargin=*]
\item\label{item:rank} If $G$ is nonregular, then $r(G)\leq3n/4$.
\item\label{item:generators} If $G$ is minimally transitive, then $d(G)\leq\log(n)$.
\item\label{item:normalizer} If $n\geq2$, then $\log([\mathbf{N}_{\Sym(n)}(G):G])\leq2n/\sqrt{\log(n)}+\log^2(n)$.
\end{enumerate}
\end{lemma}

Let $G\leq\Sym(\Omega)$ and $H\leq\Sym(\Delta)$.
A \emph{permutational isomorphism} from $G$ to $H$ consists of a group isomorphism $\varphi\colon G\to H$ and a bijection $f\colon\Omega\to\Delta$ such that for all $\omega\in\Omega$ and $g\in G$,
\[
(\omega^g)^f=(\omega^f)^{g^\varphi}.
\]
We say that $G$ and $H$ are \emph{permutationally isomorphic} if such a permutational isomorphism exists.

The following result is~\cite[Lemma~5.2]{MN}.

\begin{lemma}\label{lem:general-generators}
Let $G$ be a permutation group of degree $n\geq1$. Then $d(G)\leq n-1$.
\end{lemma}

For integers $m\geq2$ and $\ell\geq2$, let $\lpp(m)$ be the largest prime-power divisor of $m$, and let
\begin{equation}\label{eq:theta}
\theta(m,\ell)=\frac{2m\ell\log(2\ell)}{\lpp(m)}+\frac{16m\ell}{\sqrt{\log(\lpp(m))}}.
\end{equation}

\begin{lemma}\label{lem:uniform-theta}
There exists an absolute constant $C>0$ such that, for every sufficiently large integer $n$ and all integers $m\geq2$ and $\ell\geq2$ satisfying $m\ell=n$ and $\ell\leq\log^{24}(n)$,
\[
\theta(m,\ell)\leq\frac{Cn}{\sqrt{\log(\log(n))}}.
\]
\end{lemma}

\begin{proof}
Since every prime divisor of $m$ is at most $\lpp(m)$, there are at most $\lpp(m)$ distinct prime divisors of $m$, each contributing a prime-power part at most $\lpp(m)$. Therefore,
\[
m\leq2^{\lpp(m)\log(\lpp(m))}\leq2^{\lpp(m)^2}.
\]
Note $m\geq\sqrt n$ for sufficiently large $n$. Hence $\lpp(m)\geq\sqrt{\log(n)}/2$ and $\log(\lpp(m))\geq\log(\log(n))/4$.
Substituting these bounds together with $\log(2\ell)\leq24\log(\log(n))+1$ into the definition of $\theta(m,\ell)$, we obtain for some absolute constant $C>0$ that
\[
\theta(m,\ell)\leq\frac{4n(24\log(\log(n))+1)}{\sqrt{\log(n)}}+\frac{32n}{\sqrt{\log(\log(n))}}\leq\frac{Cn}{\sqrt{\log(\log(n))}}.\qedhere
\]
\end{proof}

We will also use the following bound for trivial sections of induced modules.

\begin{lemma}\label{lem:trivial-sections}
Let $H$ be a subgroup of a finite group $G$ with $m=[G:H]\geq2$, and let $A$ be a finite abelian $H$-module.
Suppose that $Y\leq X$ are $\mathbb{Z} G$-submodules of $\Ind_H^G(A)$, and $G$ acts trivially on $X/Y$. Then
\[
\log(|X/Y|)\leq\frac{2m\log(|A|)}{\sqrt{\pi\log(m)}}.
\]
\end{lemma}

\begin{proof}
First assume that $A$ is elementary abelian of order $p^\ell$, where $p$ is a prime.
By~\cite[Lemma~2.1 and Remark~2.2]{LSS2021}, $X$ is generated as an $\mathbb{F}_pG$-module by at most $2m\ell/\sqrt{\pi\log(m)}$ elements.
Since $G$ acts trivially on $X/Y$, the images of these generators span $X/Y$ over $\mathbb{F}_p$.
Thus $\log(|X/Y|)\leq2m\ell\log(p)/\sqrt{\pi\log(m)}=2m\log(|A|)/\sqrt{\pi\log(m)}$.

For general $A$, choose an $H$-invariant series $0=A_0\leq A_1\leq\cdots\leq A_t=A$ whose factors are elementary abelian.
Let $U_i=\Ind_H^G(A_i)$ for $i\in\{0,\ldots,t\}$. Then we have the series
\[
0=U_0\leq U_1\leq\cdots\leq U_t=\Ind_H^G(A)
\]
with $U_i/U_{i-1}\cong\Ind_H^G(A_i/A_{i-1})$ for $i\in\{1,\ldots,t\}$, and the submodules $B_i=((X\cap U_i)+Y)/Y$ form a series $0=B_0\leq B_1\leq\cdots\leq B_t=X/Y$.
For each $i\in\{1,\ldots,t\}$, the factor $B_i/B_{i-1}$ is a quotient of $(X\cap U_i)/(X\cap U_{i-1})$, which is isomorphic to a submodule of $U_i/U_{i-1}$.
Since $G$ acts trivially on $X/Y$, it follows that $B_i/B_{i-1}$ is a trivial section of $\Ind_H^G(A_i/A_{i-1})$.
Then applying the elementary abelian case to these factors gives
\[
\log(|X/Y|)\leq\sum_{i=1}^t\frac{2m\log(|A_i/A_{i-1}|)}{\sqrt{\pi\log(m)}}=\frac{2m\log(|A|)}{\sqrt{\pi\log(m)}}.\qedhere
\]
\end{proof}

\section{Auxiliary estimates for counting orbital partitions}\label{sec:auxiliary}

The bounds in this section will be used in the proof of Proposition~\ref{prop:general-recovery}.

\subsection{Quotients associated with a transitive action}\label{subsec:prime-parts}

Let $L$ be a finite group, let $L_1\leq L$ be core-free with $[L:L_1]=\ell\geq2$, and let $\mathcal{O}(L,L_1)$ be the set of distinct kernels of the actions of $L$ on the nontrivial irreducible constituents of the permutation $\mathbb{C} L$-module $\mathbb{C}[L/L_1]$ afforded by right multiplication on $L/L_1$. For $K\in\mathcal{O}(L,L_1)$, let $V_K$ be the sum of all irreducible $\mathbb{C} L$-submodules of $\mathbb{C}[L/L_1]$ on which the action of $L$ has kernel $K$, and let $d_K=\dim_{\mathbb{C}}(V_K)$.
Let
\[
\Lambda=\sum_{K\in\mathcal{O}(L,L_1)}\log(|L/K|),\quad
\Lambda_0=\sum_{K\in\mathcal{O}(L,L_1)}\log(|\mathbf{Z}(L/K)|),\quad
\Lambda_p=\sum_{K\in\mathcal{O}(L,L_1)}\log(|L/K|_p)
\]
for every prime $p$.

\begin{lemma}\label{lem:quotient-points}
Let $L$ be a finite group, let $L_1\leq L$ be core-free with $[L:L_1]=\ell\geq2$, and let $K\in\mathcal{O}(L,L_1)$. Then $\core_L(L_1K)=K$, and $\mathbf{Z}(L/K)$ is cyclic of order at most $\ell$.
\end{lemma}

\begin{proof}
Choose an irreducible submodule $W$ of $\mathbb{C}[L/L_1]$ with kernel $K$. Orthogonal projection onto $W$, with respect to the standard inner product, commutes with the action of $L$. Hence the projection of the basis vector corresponding to $L_1$ is $L_1$-fixed and nonzero: otherwise, transitivity would force every basis vector to have zero projection.
In the resulting faithful irreducible representation of $L/K$, the subgroup $L_1K/K$ fixes a nonzero vector.
The fixed space of $\core_{L/K}(L_1K/K)$ is nonzero and $(L/K)$-invariant. Irreducibility makes this the whole representation space, and faithfulness then gives $\core_{L/K}(L_1K/K)=1$. Thus $\core_L(L_1K)=K$.

By Schur's lemma, $\mathbf{Z}(L/K)$ acts faithfully by scalars on $W$, so $\mathbf{Z}(L/K)$ is cyclic.
The action of $L/K$ on the right cosets of $L_1K/K$ is faithful and transitive. Its centre is semiregular, so $|\mathbf{Z}(L/K)|\leq[L:L_1K]\leq\ell$.
\end{proof}

\begin{lemma}\label{lem:prime-energy}
Let $L$ be a finite group, and let $L_1\leq L$ be core-free with $[L:L_1]=\ell\geq2$. Then
\[
\Lambda\leq2(\ell-1)\log(2\ell),\quad
\Lambda_0\leq2(\ell-1),\quad \Lambda_p\leq2(\ell-1)
\]
for every prime $p$.
\end{lemma}

\begin{proof}
For each $K\in\mathcal{O}(L,L_1)$, Galois automorphisms of a cyclotomic splitting field permute the irreducible characters with kernel $K$. The sum of the corresponding central primitive idempotents is therefore rational.
Applying this idempotent to $\mathbb{Q}[L/L_1]$ gives a faithful rational representation of $L/K$ of dimension $d_K$.
Transitivity gives one trivial constituent, so $\sum_K d_K=\ell-1$.
Minkowski's theorem, in the form of~\cite[Theorem~1]{Serre}, gives the $p$-adic valuation
\[
v_p(|L/K|)\leq\sum_{j\geq0}\left\lfloor\frac{d_K}{p^j(p-1)}\right\rfloor
\leq\frac{p d_K}{(p-1)^2}.
\]
Multiplying by $\log(p)$ and summing over $K\in\mathcal{O}(L,L_1)$ gives $\Lambda_p\leq p(\ell-1)\log(p)/(p-1)^2\leq2(\ell-1)$. Moreover,
\[
\left\lfloor\frac{e}{p^j(p-1)}\right\rfloor
\leq\left\lfloor\frac{2e}{p^{j+1}}\right\rfloor.
\]
Thus the order of a finite subgroup of $\GL_e(\mathbb{Q})$ divides $(2e)!$, by Legendre's formula.
Consequently,
\[
\Lambda\leq\sum_K\log((2d_K)!)\leq\sum_K 2d_K\log(2d_K)\leq2(\ell-1)\log(2\ell).
\]
Finally, let $a=|\mathbf{Z}(L/K)|$. By Lemma~\ref{lem:quotient-points}, the group $\mathbf{Z}(L/K)$ is cyclic, and its faithful scalar action on each irreducible constituent of $V_K$ gives $\phi(a)\leq d_K$ by Galois conjugacy.
The elementary inequalities $\phi(a)\geq\sqrt{a/2}$ and $\log(a)\leq\sqrt{2a}$ give $\log(a)\leq2d_K$; the first inequality follows by multiplying the prime-power factors of $\phi(a)^2/a$.
Summing over $K\in\mathcal{O}(L,L_1)$ gives $\Lambda_0\leq2(\ell-1)$.
\end{proof}

\begin{lemma}\label{lem:joint-faithful}
Let $L$ be a finite group, let $L_1\leq L$ be core-free with $[L:L_1]=\ell\geq1$, let $x\geq0$ be a real number, let $\mathcal K$ be a set of normal subgroups of $L$ such that $\core_L(L_1K)=K$ and $\log(|L/K|)\leq x$ for every $K\in\mathcal K$, and let $N=\bigcap_{K\in\mathcal K}K$, with $N=L$ if $\mathcal K$ is empty. Then
\[
\log(|L/N|)\leq x\log(\ell).
\]
\end{lemma}

\begin{proof}
If $\mathcal K$ is empty, then $N=L$ and the conclusion holds. Assume that $\mathcal K$ is nonempty.
Since $L_1N\leq L_1K$ for every $K$, normality gives
\[
N\leq\core_L(L_1N)\leq\bigcap_K\core_L(L_1K)=N.
\]
Pass to $L/N$ and its faithful action on $L/L_1N$, whose degree is at most $\ell$.
We may therefore assume that $N=1$.
Choose an irredundant subfamily $K_1,\ldots,K_t$ with intersection one, and let
$A_i=\bigcap_{j\neq i}K_j$.
Then $A_i\neq1$ and $A_i\cap K_i=1$.
For distinct $i$ and $j$, we have $[A_i,A_j]\leq A_i\cap A_j=1$. Moreover, $A_i\cap\prod_{j\neq i}A_j\leq A_i\cap K_i=1$, so $\prod_i A_i$ is an internal direct product.
Let $V=L_1\cap\prod_i A_i$, and let $V_i$ be its projection to $A_i$.
If $V_i=A_i$, then $A_i\leq L_1K_i$, because the other factors lie in $K_i$.
Normality would imply $A_i\leq\core_L(L_1K_i)=K_i$, a contradiction.
Therefore,
\[
2^t\leq\prod_i[A_i:V_i]
\leq\bigl[\prod_i A_i:V\bigr]\leq[L:L_1]\leq \ell.
\]
Returning to the original group $L$, we have selected at most $\log(\ell)$ members $K_1,\ldots,K_t$ of $\mathcal K$ with intersection $N$. The quotient homomorphisms give a monomorphism $L/N\to\prod_{i=1}^t L/K_i$, and hence $\log(|L/N|)\leq tx\leq x\log(\ell)$.
\end{proof}

\subsection{Extensions with a prescribed restriction}\label{subsec:extensions}

We count extensions up to equivalence, so an isomorphism of extensions must induce the identity on both the kernel and the quotient. We first consider extensions with abelian kernel that split over a given subgroup.

\begin{samepage}
\begin{lemma}\label{lem:restriction}
Let $p$ be a prime, let $H$ be a subgroup of a finite group $G$, and let $A$ be a finite abelian $p$-group with a prescribed action of $G$ by automorphisms.
There are at most $|G|_p^{d(A)v_p([G:H])}$ equivalence classes of extensions
\[
1\longrightarrow A\longrightarrow E\xrightarrow{\varepsilon}G\longrightarrow1
\]
which induce this action and in which $A$ has a complement in $\varepsilon^{-1}(H)$.
\end{lemma}
\end{samepage}

\begin{proof}
The extension classes inducing the given action correspond to the elements of $\mathrm{H}^2(G,A)$, and the extension splits over $H$ precisely when its class restricts to zero in $\mathrm{H}^2(H,A)$~\cite[Chapter~IV, Theorem~3.12]{Brown}. We thus bound the kernel of the restriction homomorphism $\Res_H^G$. Set $e=d(A)$. Every elementary abelian section of a finite abelian $p$-group $W$ has order at most $p^{d(W)}$, since every subgroup of $W$ can be generated by at most $d(W)$ elements.

Let $U$ be a finite $p$-group acting on $A$, and let $V$ be a normal subgroup of $U$ of index $p$. The exact sequence of Dekimpe, Hartl and Wauters~\cite[Lemma~6.8]{DHW} contains
\[
\mathrm{H}^2(U/V,A^V)\longrightarrow\Ker(\Res_V^U)\longrightarrow\mathrm{H}^1(U/V,\mathrm{H}^1(V,A)),
\]
where restriction is taken in degree two. Since $U/V$ is cyclic of order $p$, the first group is a section of $A^V$ killed by $p$ and has order at most $p^e$, and the last group is a section of $\mathrm{H}^1(V,A)$ killed by $p$. As a one-cocycle is determined by its values on a generating set of $V$, the group $\mathrm{H}^1(V,A)$ is a quotient of a subgroup of $A^{d(V)}$, so the last group has order at most $p^{ed(V)}\leq|V|^e$. Exactness gives $|\Ker(\Res_V^U)|\leq p^e|V|^e=|U|^e$.

Choose a Sylow $p$-subgroup $P$ of $G$ containing a Sylow $p$-subgroup of $H$, and set $t=v_p([G:H])$. Then $[P:P\cap H]=p^t$. Every subgroup of a $p$-group has a chain to the group with successive normal inclusions of index $p$.
Apply the normalizer condition at each step to obtain such a chain from $P\cap H$ to $P$.
Its length is $t$.
Multiplying the restriction-kernel bounds along the chain gives at most $|P|^{et}$.

Finally, $\Res_P^G\colon\mathrm{H}^2(G,A)\to\mathrm{H}^2(P,A)$ is a monomorphism.
Restriction followed by transfer is multiplication by $[G:P]$, which is prime to $p$~\cite[Chapter~III, Section~9]{Brown}.
Naturality of restriction embeds $\Ker(\Res_H^G)$ in $\Ker(\Res_{P\cap H}^P)$, whose order is at most $|P|^{et}=|G|_p^{d(A)v_p([G:H])}$.
\end{proof}

Fix finite groups $K$ and $G$, a subgroup $H\leq G$, and an exact sequence
\[
1\longrightarrow K\xrightarrow{\alpha}V\xrightarrow{\eta}H\longrightarrow1.
\]
An \emph{extension of $G$ by $K$ with prescribed restriction $V$} is a tuple $(W,\beta,\varepsilon,\iota)$ such that
\[
1\longrightarrow K\xrightarrow{\beta}W\xrightarrow{\varepsilon}G\longrightarrow1
\]
is exact and $\iota\colon V\to W$ is a monomorphism satisfying $\alpha\iota=\beta$ and $\iota\varepsilon=\eta$, where $\eta$ is followed by the inclusion $H\hookrightarrow G$.
Two such tuples $(W_1,\beta_1,\varepsilon_1,\iota_1)$ and $(W_2,\beta_2,\varepsilon_2,\iota_2)$ are called \emph{equivalent} if there is an isomorphism $\varphi\colon W_1\to W_2$ satisfying
\[
\beta_1\varphi=\beta_2,\quad \varphi\varepsilon_2=\varepsilon_1,\quad \iota_1\varphi=\iota_2.
\]

\begin{lemma}\label{lem:relative-extensions}
Let $\mathcal E$ be the set of equivalence classes of extensions of $G$ by $K$ with prescribed restriction $V$ as defined above, let $x=\log(|G|)$, and let $y=\log(|K|)$. Then
\[
|\mathcal E|\leq2^{xy^2+xy(\log([G:H])+1)}.
\]
\end{lemma}

\begin{proof}
Choose a homomorphism from $G$ to $\Out(K)$ extending the outer action prescribed by $V$.
There are at most $|\Out(K)|^{d(G)}$ choices.
Fix one such outer action and one extension inducing it whose restriction to $H$ is equivalent to $V$. By the classification of extensions with a given outer action~\cite[Chapter~IV, Theorem~6.6]{Brown}, the equivalence classes of extensions with this outer action are parametrized by $\mathrm{H}^2(G,\mathbf{Z}(K))$, with the chosen extension corresponding to zero. This parametrization commutes with restriction to $H$. Hence the number of extension classes whose restrictions to $H$ are equivalent to $V$ is the number of equivalence classes of extensions of $G$ by $\mathbf{Z}(K)$ inducing the given action on the centre and splitting over $H$. The cohomology of the centre is the direct sum over its primary components, so Lemma~\ref{lem:restriction} bounds the logarithm of this number by
\[
x\sum_p d(\mathbf{Z}(K)_p)v_p([G:H])\leq xy\log([G:H]).
\]

For each such extension $W$, it remains to choose the monomorphism $\iota\colon V\to W$. It is determined by its values on lifts of $d(H)$ generators of $H$, since its restriction to $K$ is prescribed. For each lift, any two possible images lie in the same coset of $K$ and induce the same conjugation on $K$, so they differ by an element of $\mathbf{Z}(K)$. There are thus at most $|\mathbf{Z}(K)|^{d(H)}\leq 2^{xy}$ choices of $\iota$. Since $\log(|\Out(K)|)\leq y^2$ and $d(G)\leq x$, multiplying the bounds for the outer actions, the extension classes and the monomorphisms gives the stated inequality.
\end{proof}

\subsection{Invariant orbits of affine actions}\label{subsec:affine}

An action of a group $G$ on the underlying set of a group $K$ is called \emph{affine} if there is a homomorphism $\alpha\colon G\to\Aut(K)$ such that
\[
v^g=1_K^g\,v^{\alpha(g)}
\]
for every $g\in G$ and $v\in K$. The homomorphism $\alpha$ is uniquely determined by this affine action and is called its \emph{associated homomorphism}. For each $g\in G$, we call $\alpha(g)$ its \emph{linear part} and the permutation $v\mapsto1_K^g v$ its \emph{translation part}.

\begin{lemma}\label{lem:central-fibres}
Let $K$ be a finite nilpotent group of class $c$, let $G$ act affinely on $K$ with associated homomorphism $\alpha$, let $G_0\trianglelefteq G$ with $\alpha(G_0)\leq\Inn(K)$, and let $M=G/G_0$.
Then $\alpha$ induces an action of $M$ on each factor of the upper central series $1=\mathbf{Z}_0(K)\leq\cdots\leq\mathbf{Z}_c(K)=K$.
For each $i\in\{0,\ldots,c-1\}$, let $b_i\geq1$ be a real number such that every trivial $M$-section of $\mathbf{Z}_{i+1}(K)/\mathbf{Z}_i(K)$ has order at most $b_i$.
Then the number of $G_0$-orbits on $K$ that are invariant under $G$ is at most $\prod_{i=0}^{c-1}b_i$.
\end{lemma}

\begin{proof}
The assertion is immediate if $K=1$. Assume that $K\neq1$. Since inner automorphisms of $K$ act trivially on each upper central factor, the action of $G$ on each such factor induced by $\alpha$ factors through $M$.
It suffices to bound the number of $G$-invariant $G_0$-orbits above a given $G$-invariant $G_0$-orbit on $K/\mathbf{Z}(K)$, and then to induct on the nilpotency class.
Choose $\overline{v}\in K/\mathbf{Z}(K)$ in that orbit.
Since $G$ stabilizes this orbit, each lift of an element of $M$ can be adjusted by an element of $G_0$ to fix $\overline{v}$.
Thus the stabilizer $G_{\overline{v}}$ of $\overline{v}$ maps onto $M$, with kernel $(G_0)_{\overline{v}}$.

Choose $v\in K$ such that $\overline{v}=v\mathbf{Z}(K)$, and identify the fibre $v\mathbf{Z}(K)$ with $\mathbf{Z}(K)$ by $vz\mapsto z$.
The group $(G_0)_{\overline{v}}$ acts on this fibre by translations, since $\alpha(g)$ fixes $\mathbf{Z}(K)$ pointwise for every $g\in(G_0)_{\overline{v}}$.
Let $J=\{v^{-1}v^g\mid g\in(G_0)_{\overline{v}}\}\leq \mathbf{Z}(K)$. This subgroup is invariant under the action of $G_{\overline{v}}$ induced by $\alpha$.
The $(G_0)_{\overline{v}}$-orbits in the fibre are therefore indexed by $\mathbf{Z}(K)/J$, and correspond bijectively to the $G_0$-orbits above the chosen orbit on $K/\mathbf{Z}(K)$.
The induced action of $G_{\overline{v}}$ on $\mathbf{Z}(K)/J$ is affine and factors through $M$.
If its fixed-point set is nonempty, then translating one fixed point to the identity identifies that set with $(\mathbf{Z}(K)/J)^M$.
This is a trivial $M$-section of $\mathbf{Z}(K)$, so there are at most $b_0$ invariant orbits above the chosen orbit on $K/\mathbf{Z}(K)$.
Induction applied to $K/\mathbf{Z}(K)$ gives at most $\prod_{i=1}^{c-1}b_i$ choices for that orbit, and hence the required bound.
\end{proof}

\begin{lemma}\label{lem:refined-affine}
Let $M$ be a finite group with a subgroup $M_1$ of index $m\geq2$, let $\varepsilon\colon G\to M$ be a surjective homomorphism with kernel $G_0$, and let $t\geq1$ be an integer. For each $i\in\{1,\ldots,t\}$, let $M_{1,i}$ be a subgroup of $M_1$ of index $m_i$, and let $B_i$ be a finite group. Let
\[
B=\prod_{i=1}^t \prod_{j\in M/M_{1,i}} B_i,\quad
\Psi=\sum_{i=1}^t m_i\log(|B_i|),\quad
\Psi_p=\sum_{i=1}^t m_i\log(|B_i|_p)
\]
for every prime $p$.
Suppose that $G$ acts affinely on $B$ with associated homomorphism $\alpha$ satisfying $\alpha(G_0)\leq\Inn(B)$, and each $\alpha(g)$ with $g\in G$ maps the coordinate factor indexed by $(i,j)$ isomorphically onto the factor indexed by $(i,j^{\varepsilon(g)})$ for all $i\in\{1,\ldots,t\}$ and $j\in M/M_{1,i}$. Then for every prime divisor $p$ of $m$, the number of $G_0$-orbits on $B$ invariant under $G$ is at most
\[
2^{\frac{m\Psi}{m_p}+\frac{8m\Psi_p}{\sqrt{\log(m_p)}}}.
\]
\end{lemma}

\begin{proof}
Choose a Sylow $p$-subgroup $T$ of $M$ and a Sylow $p$-subgroup $P$ of its preimage in $G$.
Let $P_0=P\cap G_0$.
Then $P/P_0\cong T$, and we identify this quotient with $T$.
For each $i\in\{1,\ldots,t\}$, the permutations $v\mapsto a^{-1}vb$ of $B_i$, with $a\in B_i$ and $b\in B_i$, form a group isomorphic to
\[
D_i=(B_i\times B_i)/\{(z,z)\mid z\in \mathbf{Z}(B_i)\}.
\]
Let $D=\prod_{i=1}^t \prod_{j\in M/M_{1,i}} D_i$, acting coordinatewise on $B$. The affine action of $G$ normalizes $D$, and $\alpha(G_0)\leq\Inn(B)$ implies that the image of $G_0$ lies in $D$.
The $p$-group $P$ normalizes a Sylow $p$-subgroup $K$ of $D$, because the number of Sylow $p$-subgroups of $D$ is prime to $p$.
The image of $P_0$ lies in $K$.
Indeed, it is a $p$-subgroup of $D$ normalizing $K$.
For each $i\in\{1,\ldots,t\}$ and $j\in M/M_{1,i}$, let $K_{i,j}$ be the projection of $K$ to the copy of $D_i$ indexed by $(i,j)$. Then $K$ is the direct product of these Sylow $p$-subgroups, and
\begin{equation}\label{eq:coordinate-sylow}
\log(|K_{i,j}|)=2\log(|B_i|_p)-\log(|\mathbf{Z}(B_i)|_p)\leq2\log(|B_i|_p).
\end{equation}

Let $A$ be the preimage of $T$ in $G$.
Every $G_0$-orbit stabilized by $A$ contains a $P_0$-orbit stabilized by $P$.
To see this, choose a point $v$ of such a $G_0$-orbit.
Since $G_0$ is transitive on this orbit and $A$ stabilizes it, $A=G_0 A_v$.
A Sylow $p$-subgroup $S$ of $A_v$ maps onto $T$.
It is contained in a Sylow $p$-subgroup of $A$, and the Sylow subgroups of $A$ are conjugate by $G_0$: indeed $A=G_0 P$, so a conjugating element can be modified by an element of $P$ to lie in $G_0$.
After this conjugation, $S\leq P$ fixes the resulting point $v'$ and $P=P_0S$.
Hence the orbit $(v')^{P_0}$ is $P$-invariant.
Distinct $G_0$-orbits yield distinct chosen $P_0$-orbits.

Each such orbit lies in a $P$-invariant $K$-orbit.
For $j\in M/M_{1,i}$, we have $|T_j|\leq |M_{1,i}|_p$, and hence $|j^T|\geq[M:M_{1,i}]_p\geq m_p$.
Since $P$ permutes the coordinate $K_{i,j}$-orbits, a $P$-invariant $K$-orbit is determined by its coordinate orbit at one representative of each $T$-orbit $O$.
There are at most $|B_i|$ choices for each $O$ of type $i$.
Hence there are at most $2^{m\Psi/m_p}$ such $P$-invariant $K$-orbits.

Fix such a $K$-orbit $X$ and $x\in X$.
Let $\phi_s$ be the automorphism $k\mapsto s^{-1}ks$ of $K$, for $s\in P$, and let
\[
\widehat P=\{(s,c)\mid s\in P,\ c\in K,\ x^s=x^c\}.
\]
With multiplication $(s,c_1)(t,c_2)=(st,c_2c_1^{\phi_t})$ this group acts on $K$ by $k^{(s,c)}=ck^{\phi_s}$.
Let $\widehat P_0$ be the preimage of $P_0$ under the projection $\widehat P\to P$ to the first coordinate. The linear part of each element of $\widehat P_0$ is an inner automorphism of $K$, because the image of $P_0$ lies in $K$.
The mapping $k\mapsto x^k$ induces a $T$-equivariant bijection between the set of $\widehat P_0$-orbits on $K$ and the set of $P_0$-orbits on $X$.
Apply Lemma~\ref{lem:central-fibres} to the affine action of $\widehat P$ on $K$, with normal subgroup $\widehat P_0$ and associated homomorphism $(s,c)\mapsto\phi_s$.
For each $i\in\{1,\ldots,t\}$ and each $T$-orbit $O$ on $M/M_{1,i}$, choose $j\in O$. Each upper central factor of $\prod_{h\in O}K_{i,h}$ is induced from the corresponding factor of $K_{i,j}$, regarded as a $T_j$-module. The index $[T:T_j]$ is $|O|$.
By Lemma~\ref{lem:trivial-sections}, the product of the bounds on the orders of their trivial $T$-sections is at most
\[
2^{\frac{2|O|\log(|K_{i,j}|)}{\sqrt{\pi\log(|O|)}}}.
\]
Multiplying these bounds over the coordinate orbits and using~\eqref{eq:coordinate-sylow} and $|O|\geq m_p$, we obtain at most $2^{4m\Psi_p/\sqrt{\pi\log(m_p)}}\leq2^{8m\Psi_p/\sqrt{\log(m_p)}}$ $P_0$-orbits in $X$ stabilized by $P$.
Multiplying by the bound $2^{m\Psi/m_p}$ for the number of $P$-invariant $K$-orbits proves the lemma.
\end{proof}

\section{Orbital partitions with a prescribed quotient action}\label{sec:recovery}

A key step in the counting argument is to bound the number of orbital partitions when the action on the orbits of a normal subgroup and the possible actions of that subgroup on one orbit are prescribed.
We establish bounds of this kind for use in Section~\ref{sec:counting}.

Throughout this section, fix integers $m\geq2$ and $\ell\geq2$, a transitive permutation group $M$ on $\{1,\ldots,m\}$, with point stabilizer $M_1$ of $1$, and a nonempty set $\mathcal L$ of transitive permutation groups on $\{1,\ldots,\ell\}$, no two of which are permutationally isomorphic.
Consider a transitive permutation group $G$ on a set $\Omega$ of size $m\ell$ and a normal subgroup $K$ of $G$.
Suppose that $K$ is the kernel of the action of $G$ on the set of $K$-orbits and that these orbits can be labelled by $M/M_1$ so that the induced permutation group is $M$ in its action on $M/M_1$. Under this labelling we identify $G/K$ with $M$ and write $\overline{g}$ for the image of $g\in G$ in $M$. Let $\sigma\colon\Omega\to M/M_1$ send each point to the label of its $K$-orbit. Then for all $\omega\in\Omega$ and $g\in G$,
\begin{equation}\label{eq:quotient-action}
(\omega^g)^\sigma=(\omega^\sigma)^{\overline{g}}.
\end{equation}
We require the permutation group induced by $K$ on $\sigma^{-1}(M_1)$ to be permutationally isomorphic to a member $L\in\mathcal L$.

As $G$, $K$, $\Omega$ and $\sigma$ vary subject to these conditions, we count the pairs $(\mathcal P,\sigma)$, where $\mathcal P$ is the orbital partition of $G$, up to isomorphisms preserving the labelling of the $K$-orbits by $M/M_1$.
For this purpose, we call two such pairs $(\mathcal P_1,\sigma_1)$ and $(\mathcal P_2,\sigma_2)$ on point sets $\Omega_1$ and $\Omega_2$, respectively, \emph{isomorphic} if there is an isomorphism $f\colon\Omega_1\to\Omega_2$ from $\mathcal P_1$ to $\mathcal P_2$ in the sense of Definition~\ref{def:orbital-isomorphism} such that $f\sigma_2=\sigma_1$.
Thus $f$ carries $\sigma_1^{-1}(v)$ onto $\sigma_2^{-1}(v)$ for every $v\in M/M_1$.
Let $\mathcal R=\mathcal R(M,\mathcal L)$ be the set of isomorphism classes of the pairs $(\mathcal P,\sigma)$ arising from some choice of $G$, $K$, $\Omega$ and $\sigma$ satisfying the above conditions.
The set $\mathcal R$ is nonempty, as for $L\in\mathcal L$ we may take $G=L\wr M$ in its natural imprimitive action with base group $K=L^m$.

Including $\sigma$ fixes the labelling of the $K$-orbits by $M/M_1$, allowing the constructions and counts below to use coordinates indexed by this fixed point set.
Passing from $(\mathcal P,\sigma)$ to $\mathcal P$ induces a surjection from $\mathcal R$ onto the set of isomorphism classes of orbital partitions arising from the groups described above.
Hence the bounds for $|\mathcal R|$ below also bound the number of these classes.

\begin{proposition}\label{prop:general-recovery}
Let $m,\ell,M,\mathcal L$ and $\mathcal R=\mathcal R(M,\mathcal L)$ be as above, and let $n=m\ell$.
There exists an absolute constant $C$ such that the following inequalities hold:
\begin{align}
\log(|\mathcal R|)\leq{}&C\ell\log(\ell)\log(|M|)\bigl(\ell\log(\ell)+\log(n)+1\bigr)\notag\\
&+C\ell\log^2(n)+\theta(m,\ell)+\frac{8n}{\sqrt{\log(m)}}.\label{eq:general-full}
\end{align}
\begin{align}
\log(|\mathcal R|)\leq{}&C\bigl(\ell^2\log^2(\ell)\log(n)+\ell\log^2(n)\bigr)\notag\\
&+d(M_1)\max_{L\in\mathcal L}\log([\mathbf{N}_{\Sym(\ell)}(L):L])+\theta(m,\ell)+Cd(M)m\sqrt \ell\,\log^{\frac{3}{2}}(\ell).\label{eq:general-outer}
\end{align}
\end{proposition}

In the remainder of this section we develop the proof of Proposition~\ref{prop:general-recovery} in several steps, bringing the estimates together at the end of the section to complete the argument.

\subsection{The action of the block kernel on ordered pairs}

Fix $G,K,\Omega$ and $\sigma$ as in Proposition~\ref{prop:general-recovery}. Let $\overline{\Omega}$ be the set of $K$-orbits on $\Omega$, identified with $M/M_1$ via the labelling given by $\sigma$, and write $\Omega_v=\sigma^{-1}(v)$ for $v\in\overline{\Omega}$. Let $L_v$ be the group induced by $K$ on $\Omega_v$, and let $\zeta_v\colon K\to L_v$ be the induced homomorphism. Thus $K$ is a subdirect subgroup of $\prod_{v\in\overline{\Omega}}L_v$. Even when all the groups $L_v$ are known, the actions of an element of $K$ on two different blocks need not be independent. We first show which relations between the block actions suffice to determine the orbits of $K$ on the set $\Omega^2$ of ordered pairs.

Consider a nontrivial irreducible constituent of $\mathbb C[\Omega_v]$ as an $L_v$-module, and let $R$ be its kernel in $L_v$. As a $K$-module, it has kernel $D=\zeta_v^{-1}(R)$. The value of $k\in K$ in $K/D$ can therefore be read from its action on $\Omega_v$ by the epimorphism
\[
L_v\longrightarrow K/D,\quad \zeta_v(k)\longmapsto kD.
\]
If the same $D$ occurs on another block, then the two block actions of $k$ must give the same value in $K/D$. Now let $v$ range over $\overline{\Omega}$ and, for each $v$, let the constituent range over all nontrivial irreducible constituents of $\mathbb C[\Omega_v]$. For every kernel $D$ obtained in this way, impose these equalities on the coordinates of $\prod_{v\in \overline{\Omega}}L_v$, and let $S$ be the subgroup they define. Thus $s=(s_v)_{v\in \overline{\Omega}}$ lies in $S$ precisely when, for every such $D$, the images of $s_v$ in $K/D$ agree on all blocks on which $D$ occurs.

\begin{lemma}\label{lem:block-pair-orbits}
With the notation above, $K\leq S$, and $K$ and $S$ have the same orbits on $\Omega^2$.
\end{lemma}

\begin{proof}
Every $k\in K$ satisfies the defining equalities of $S$, and so $K\leq S$. Conversely, fix $s\in S$. For each kernel $D$ occurring in the construction of $S$, let $V_D$ be the sum of the irreducible $K$-submodules of $\mathbb C[\Omega]$ with kernel $D$. Choose $k_D\in K$ representing the common value of the coordinates of $s$ in $K/D$. On each block where $D$ occurs, $s_v$ and $\zeta_v(k_D)$ agree modulo $\zeta_v(D)$, so they act identically on all the constituents with kernel $D$ on that block. Since $\mathbb C[\Omega]$ is the direct sum of the block permutation modules, $s$ and $k_D$ agree on $V_D$. Both $K$ and $S$ fix the block-constant subspace pointwise.

Let $\varphi\in\End_K(\mathbb C[\Omega])$. As described in Section~\ref{sec:prelim}, $\varphi$ preserves each $V_D$, as well as $V_K$, the subspace fixed pointwise by $K$, which is the block-constant subspace. For $x\in V_D$, since $\varphi(x)\in V_D$ and $s$ agrees with $k_D$ on $V_D$, we have
\[
\varphi(x^s)=\varphi(x^{k_D})=\varphi(x)^{k_D}=\varphi(x)^s.
\]
The same equality holds on the block-constant subspace, since $s$ fixes it pointwise. These subspaces sum to $\mathbb C[\Omega]$, so $s$ commutes with every member of $\End_K(\mathbb C[\Omega])$. By Lemma~\ref{lem:orbital-endomorphisms}, $s$ stabilizes every $K$-orbit on $\Omega^2$. Thus, as $K\leq S$, the two groups have the same orbits on $\Omega^2$.
\end{proof}

Since $K\trianglelefteq G$, the quotient $M=G/K$ acts on the set of $K$-orbits on $\Omega^2$: an element of $M$ sends a $K$-orbit to its image under any lift in $G$. Each $G$-orbit on $\Omega^2$ is the union of the $K$-orbits in one $M$-orbit. By Lemma~\ref{lem:block-pair-orbits}, these $K$-orbits are precisely the $S$-orbits. Thus the quotient homomorphisms defining $S$, together with the induced action of $M$ on the set of $S$-orbits, suffice to recover the orbital partition of $G$. To count the possible quotient homomorphisms, we first determine which kernels occur on one block and how their conjugates occur on the other blocks.

\subsection{Abstract quotient data and their number}\label{subsec:one-block-quotients}

We first define and count the data for the quotient homomorphisms in Lemma~\ref{lem:block-pair-orbits}, independently of a group $G$. The data specify the quotients on one block, the labels that will identify equal kernels on different blocks, and the compatible actions on these quotients. Their realization by the groups counted in $\mathcal R$ will be proved in Subsection~\ref{subsec:realization}.

In the abstract constructions, let $\overline{\Omega}$ be the fixed set of block labels, identified with $M/M_1$, and let $\Delta=\{1,\ldots,\ell\}$. Choose $L\in\mathcal L$, let $L_1$ be its stabilizer of $1\in\Delta$, and identify $\Delta$ with $L/L_1$. Choose a homomorphism
\[
\rho\colon M_1\longrightarrow\mathbf{N}_{\Sym(\Delta)}(L)/L,
\]
and define
\begin{equation}\label{eq:abstract-block-normalizer}
N=\{(x,y)\in M_1\times\mathbf{N}_{\Sym(\Delta)}(L)\mid\rho(x)=yL\}.
\end{equation}
The projection to the first coordinate maps $N$ onto $M_1$ with kernel $1\times L$, which we identify with $L$. The projection to the second coordinate gives the action of $N$ on $\Delta$, and $N/L=M_1$.

The set $\mathcal O(L,L_1)$ consists of the distinct kernels of the nontrivial irreducible constituents of $\mathbb C[\Delta]$, as in Subsection~\ref{subsec:prime-parts}. Conjugation by $N$ permutes these kernels, and $L$ stabilizes each of them under conjugation, giving an action of $M_1$ on $\mathcal O(L,L_1)$. For each pair $(L,\rho)$, fix an ordered list $R_1,\ldots,R_t$ of representatives of these orbits. For $j\in\{1,\ldots,t\}$, let $M_{1,j}$ be the stabilizer of $R_j$ in $M_1$, and let $m_j=[M_1:M_{1,j}]$. The preimage of $M_{1,j}$ under $N\to M_1$ is $\mathbf{N}_N(R_j)$.
Write $\Lambda$, $\Lambda_0$ and $\Lambda_p$ for the sums defined in Subsection~\ref{subsec:prime-parts} for $(L,L_1)$. Since the $m_j$ kernels in the orbit of $R_j$ all have the same index in $L$, we have
\begin{equation}\label{eq:general-incidence-energy}
\sum_{j=1}^t m_j=|\mathcal O(L,L_1)|\leq\ell-1,\quad
\sum_{j=1}^t m_j\log(|L/R_j|)=\Lambda.
\end{equation}

Choose an equivalence relation $\sim$ on $\{1,\ldots,t\}$, and let $I$ consist of the least index in each class. For each $i\in I$, set $F_i=L/R_i$ and choose a subgroup $Q_i$ of $M$ containing $M_{1,i}$. For each possible $Q_i$, fix a right transversal of $Q_i$ in $M$ containing $1$. For each $j\sim i$, choose a coset $Q_i w_j$, with $w_j$ its representative in this transversal, and an epimorphism $\varepsilon_j\colon L\to F_i$ satisfying
\[
w_j M_{1,j} w_j^{-1}\leq Q_i,\quad\Ker(\varepsilon_j)=R_j.
\]
Require $w_i=1$ and let $\varepsilon_i$ be the natural quotient homomorphism $L\to L/R_i=F_i$. Thus $\varepsilon_j$ identifies $L/R_j$ with $F_i$.

For $j\sim i$, conjugation by $\mathbf{N}_N(R_j)$ on $L/R_j$, followed by this identification, induces an outer action $\beta_j\colon M_{1,j}\to\Out(F_i)$. Explicitly, for $x\in M_{1,j}$ choose a preimage $h\in\mathbf{N}_N(R_j)$ and take the outer class of the automorphism induced by $zR_j\mapsto h^{-1}zhR_j$. This is independent of $h$, since the kernel $L$ of $\mathbf{N}_N(R_j)\to M_{1,j}$ induces inner automorphisms. We require the homomorphisms $\varphi_i\colon Q_i\to\Out(F_i)$ used below to satisfy
\begin{equation}\label{eq:compatible-outer-actions}
\varphi_i(w_j x w_j^{-1})=\beta_j(x)\quad\text{for all }j\sim i\text{ and }x\in M_{1,j}.
\end{equation}

We are now ready to define the two sets $\mathcal T_1$ and $\mathcal T_2$ of tuples used in the counting argument. Let $\mathcal T_1$ be the set of all tuples
\[
\bigl(L,\rho,\sim,(Q_i,E_i)_{i\in I},(Q_i w_j,\varepsilon_j)_{j=1}^t\bigr),
\]
where the common data range over the choices above and each $E_i$ ranges over fixed representatives of equivalence classes of extensions
\[
1\longrightarrow F_i\longrightarrow E_i\longrightarrow Q_i\longrightarrow1
\]
with the following two properties: the induced outer action $\varphi_i$ satisfies~\eqref{eq:compatible-outer-actions}, and the preimage of $M_{1,i}$ under $E_i\to Q_i$ is equivalent to $\mathbf{N}_N(R_i)/R_i$ as an extension of $M_{1,i}$ by $F_i$. Here $\varphi_i(q)$ is the outer class of conjugation on $F_i$ by any preimage of $q$ in $E_i$, and extension equivalences induce the identity on the specified kernel and quotient.
Let $\mathcal T_2$ be the set of all tuples
\[
\bigl(L,\rho,\sim,(Q_i,\varphi_i)_{i\in I},(Q_i w_j,\varepsilon_j)_{j=1}^t\bigr),
\]
where the common data range over the same choices and each $\varphi_i\colon Q_i\to\Out(F_i)$ satisfies~\eqref{eq:compatible-outer-actions}. In particular, a tuple in $\mathcal T_2$ includes no choice of an extension $E_i$. Neither definition requires a tuple to arise from a permutation group $G$.

\begin{lemma}\label{lem:recovery-tuples}
Let $\mathcal T_1$ and $\mathcal T_2$ be the sets of tuples defined above, and let $n=m\ell$.
There exists an absolute constant $C$ such that the following inequalities hold:
\begin{equation}\label{eq:general-full-front}
|\mathcal T_1|\leq
2^{C\bigl(\ell^2\log^2(\ell)\log(|M|)+\ell\log(\ell)\log(|M|)\log(n)\bigr)}.
\end{equation}
\begin{equation}\label{eq:general-outer-front}
|\mathcal T_2|\leq
2^{C\bigl(\ell^2\log^2(\ell)\log(n)+\ell\log^2(n)\bigr)+d(M_1)\max_{L\in\mathcal L}\log([\mathbf{N}_{\Sym(\ell)}(L):L])}.
\end{equation}
\end{lemma}

\begin{proof}
Fix $L\in\mathcal L$. Let $\ell_L=\log(|L|)\leq\ell\log(\ell)$, and let $\ell_{L,U}=\log(|L/U|)$ for each $U\in\mathcal{O}(L,L_1)$.
Since every subgroup of $\Sym(\ell)$ is generated by at most $\ell-1$ elements by Lemma~\ref{lem:general-generators}, the number of subgroups of $\Sym(\ell)$ is at most $(\ell!)^{\ell-1}$, and hence $|\mathcal L|\leq(\ell!)^{\ell-1}\leq2^{\ell^2\log(\ell)}$. Since $|M|\geq m\geq2$, the contribution of $\log(|\mathcal L|)$ is bounded by $\ell^2\log^2(\ell)\log(|M|)$ in~\eqref{eq:general-full-front} and by $\ell^2\log^2(\ell)\log(n)$ in~\eqref{eq:general-outer-front}.
Since $\log(|\Aut(L/U)|)\leq\ell_{L,U}^2$ for every $U\in\mathcal O(L,L_1)$, Lemma~\ref{lem:prime-energy} gives
\begin{equation}\label{eq:general-local-costs}
\sum_U\ell_{L,U}\leq4\ell\log(\ell),\quad \sum_U\log(|\Aut(L/U)|)\leq\sum_U\ell_{L,U}^2\leq\ell_L\sum_U\ell_{L,U}\leq4\ell^2\log^2(\ell).
\end{equation}

There are at most $[\mathbf{N}_{\Sym(\Delta)}(L):L]^{d(M_1)}$ choices of $\rho$.
For $\mathcal T_1$, we use $d(M_1)\leq\log(|M|)$ to bound this by $2^{\ell\log(\ell)\log(|M|)}$.
For $\mathcal T_2$, we retain $[\mathbf{N}_{\Sym(\Delta)}(L):L]^{d(M_1)}$.
There are at most $\ell^\ell$ partitions of $\{1,\ldots,t\}$.
For each $i$, the subgroup $Q_i$ contains $M_{1,i}$, and $[M:M_{1,i}]=mm_i\leq n$.
Hence $Q_i$ is generated by $M_{1,i}$ and at most $\log(n)$ additional elements.
These additional elements may be chosen from a fixed right transversal of $M_{1,i}$ in $M$.
Consequently, there are at most $2^{|I|\log^2(n)}\leq2^{\ell\log^2(n)}$ choices for $(Q_i)_{i\in I}$, as $|I|\leq t\leq\sum_{j=1}^t m_j\leq\ell-1$ by~\eqref{eq:general-incidence-energy}.
Moreover, there are at most $[M:Q_i]\leq n$ choices for each coset $Q_i w_j$.
Since each $\varepsilon_j$ induces an isomorphism from $L/R_j$ to $F_i$, the homomorphisms $(\varepsilon_j)_{j=1}^t$ have at most $\prod_{j=1}^t|\Aut(L/R_j)|$ choices.

For $\mathcal T_1$, the preimage of $M_{1,i}$ under $E_i\to Q_i$ and $\mathbf{N}_N(R_i)/R_i$ must define equivalent extensions of $M_{1,i}$ by $F_i$. Each such $E_i$ therefore admits an embedding of $\mathbf{N}_N(R_i)/R_i$ respecting the extension homomorphisms.
Thus, forgetting the embedding in Lemma~\ref{lem:relative-extensions} gives at most
\[
2^{\log(|M|)\ell_{L,R_i}^2+\log(|M|)\ell_{L,R_i}(\log(n)+1)}
\]
choices for $E_i$. Multiplying these bounds over $i\in I$ and including the choices of $L$, $\rho$, $\sim$, $Q_i$, $Q_i w_j$ and $\varepsilon_j$ gives~\eqref{eq:general-full-front} by~\eqref{eq:general-local-costs}, using $\log(n)=\log(m)+\log(\ell)\leq2\log(\ell)\log(|M|)$.

For $\mathcal T_2$, equation~\eqref{eq:compatible-outer-actions} in particular prescribes the restriction of $\varphi_i$ to $M_{1,i}$.
Since $Q_i$ is generated over $M_{1,i}$ by at most $\log(n)$ elements, there are at most $|\Out(F_i)|^{\log(n)}$ choices of $\varphi_i$.
Multiplying these bounds by the bounds for the choices of $L$, $\rho$, the relation $\sim$, $(Q_i)_{i\in I}$, $(Q_i w_j)_{j=1}^t$ and $(\varepsilon_j)_{j=1}^t$, and using~\eqref{eq:general-local-costs}, gives~\eqref{eq:general-outer-front}.
\end{proof}

\subsection{Matching quotient values across the blocks}

Fix an arbitrary tuple in $\mathcal T_1$ or $\mathcal T_2$. We construct families of quotient homomorphisms from these data and bound the number of invariant orbits of families. No realization by a group $G$ is assumed. We work on $\overline{\Omega}\times\Delta$, with $(v,\delta)^\sigma=v$, and let $L_v$ be a copy of $L$ acting on $\{v\}\times\Delta$.

The indices for the homomorphisms form the disjoint union of $M/M_{1,j}$ over $j\in\{1,\ldots,t\}$. For $j\sim i$ with $i\in I$, we use $a=M_{1,j}w\in M/M_{1,j}$ as the index, with label $(v,c)\in\overline{\Omega}\times M/Q_i$ given by
\begin{equation}\label{eq:general-root-map}
(v,c)=(M_1w,Q_i w_j w).
\end{equation}
Here $v$ labels the block and $c$ labels the kernel. This label is well-defined because $M_{1,j}\leq M_1$ and $w_j M_{1,j} w_j^{-1}\leq Q_i$. The sets $M/Q_i$ are also taken in disjoint union over $i$. Equal labels $c$ specify which quotient values are to be equated. For each index $a\in M/M_{1,j}$ with label $(v,c)$, where $j\sim i$ and $i\in I$, choose an epimorphism $\mu_a\colon L_v\to F_i$. Let $\mu=(\mu_a)_a$ be the family of all these epimorphisms, and define $S(\mu)$ to consist of the elements $x=(x_v)_{v\in\overline{\Omega}}\in\prod_{v\in\overline{\Omega}}L_v$ satisfying
\begin{equation}\label{eq:quotient-value-equalities}
\mu_a(x_u)=\mu_b(x_v)\quad\text{whenever $a$ and $b$ have labels $(u,c)$ and $(v,c)$, respectively.}
\end{equation}
For $x\in S(\mu)$ and $c\in M/Q_i$, define $\xi_c(x)=\mu_a(x_v)$ using any index $a$ with label $(v,c)$; the equalities~\eqref{eq:quotient-value-equalities} make this independent of the choice of $a$.

To specify and compare the families $\mu$, we use a group $Y$ acting on both the blocks and the quotients. Choose $r_v\in M$ with $M_1r_v=v$ for each $v\in\overline{\Omega}$, taking $r_{M_1}=1$. For each $c\in M/Q_i$, let $d_c$ be its representative in the right transversal of $Q_i$ in $M$ fixed in Subsection~\ref{subsec:one-block-quotients}, so that $Q_i d_c=c$. For $\mathcal T_1$, set $B_i=E_i$ and let $\tau_i\colon B_i\to\Aut(F_i)$ be conjugation on $F_i$. For $\mathcal T_2$, construct $B_i$ from the outer action $\varphi_i$ by setting
\[
B_i=\{(q,s)\in Q_i\times\Aut(F_i)\mid s\Inn(F_i)=\varphi_i(q)\},
\]
with $\tau_i\colon B_i\to\Aut(F_i)$ the projection to the second coordinate. In both cases write $D_i$ for the kernel of the projection $B_i\to Q_i$ to the first coordinate. Then $D_i=F_i$ for $\mathcal T_1$ and $D_i=\Inn(F_i)$ for $\mathcal T_2$. Let $Y$ consist of the tuples
\[
y=\bigl(w,(n_v)_{v\in\overline{\Omega}},(q_c)_{i\in I,\ c\in M/Q_i}\bigr),
\]
where $w\in M$, the image of $n_v\in N$ under $N\to M_1$ is $r_v w r_{v^w}^{-1}$, and the image of $q_c\in B_i$ under $B_i\to Q_i$ is $d_c w d_{c^w}^{-1}$. In the product $yy'$, the coordinates are $ww'$, $n_v n'_{v^w}$ and $q_c q'_{c^w}$. These multiplication rules define a group, and the projection to the first coordinate maps $Y$ onto $M$. Its kernel is $Y_0=\bigl(\prod_{v\in\overline{\Omega}}L_v\bigr)\times\prod_{i\in I}\prod_{c\in M/Q_i}D_i$. On $\overline{\Omega}\times\Delta$, the element $y$ acts by $(v,\delta)^y=(v^w,\delta^{n_v})$. On families $\mu=(\mu_a)_a$ of epimorphisms, it acts by
\begin{equation}\label{eq:general-array-action}
(\mu^y)_{a^w}=\Inn(n_v^{-1})\mu_a\tau_i(q_c).
\end{equation}
Here $a\in M/M_{1,j}$ has label $(v,c)$, where $j\sim i$ and $i\in I$.

Write $[\mu_a]=\mu_a\Inn(F_i)$. At $a=M_{1,j}$ prescribe $[\varepsilon_j]$, and at $a=M_{1,j}w$ prescribe the image of $[\varepsilon_j]$ under any $y\in Y$ whose first coordinate is $w$, using~\eqref{eq:general-array-action}. To verify that this is well-defined, note that each $y\in Y_0$ acts on $\mu_a$ by precomposition with an inner automorphism of $L_v$ and postcomposition with an inner automorphism of $F_i$, since $y\in Y_0$ implies $(\mu^y)_a=\Inn(n_v^{-1})\mu_a\tau_i(q_c)$ and $q_c\in D_i$.
It follows that $Y_0$ fixes every class $[\mu_a]$, and hence the image of a class under $y$ depends only on the first coordinate $w$. For $w\in M_{1,j}$, condition~\eqref{eq:compatible-outer-actions} says that the induced outer automorphisms $\beta_j(w)$ and $\varphi_i(w_j w w_j^{-1})$ agree. By~\eqref{eq:general-array-action}, every such $y$ therefore preserves $[\varepsilon_j]$ at $a=M_{1,j}$. Consequently, the prescribed class at $a=M_{1,j}w$ depends only on the coset $a$, and the action of $Y$ takes prescribed classes to prescribed classes.

Let $\mathcal E$ be the set of all families $\mu=(\mu_a)_a$ obtained by choosing, independently at each index $a$, an epimorphism $\mu_a$ in the class prescribed above. The action~\eqref{eq:general-array-action} therefore preserves $\mathcal E$. For $\mu\in\mathcal E$, define
\[
S_0(\mu)=
\begin{cases}
\{(x,(\xi_c(x))_{i\in I,\ c\in M/Q_i})\mid x\in S(\mu)\}&\text{if the tuple lies in }\mathcal T_1\\
\{(x,(\Inn(\xi_c(x)))_{i\in I,\ c\in M/Q_i})\mid x\in S(\mu)\}&\text{if the tuple lies in }\mathcal T_2.
\end{cases}
\]
Projection onto the first factor of $Y_0$ maps $S_0(\mu)$ isomorphically onto $S(\mu)$.

\begin{lemma}\label{lem:compatible-quotient-maps}
Fix a tuple in $\mathcal T_1$ or $\mathcal T_2$. The number of $Y$-invariant $Y_0$-orbits on $\mathcal E$ is at most $2^{\theta(m,\ell)}$. For every $\mu\in\mathcal E$, we have $S_0(\mu)\leq(Y_0)_\mu$ and $S_0(\mu)\trianglelefteq Y_\mu$.
\end{lemma}

\begin{proof}
Choose a family $\mu=(\mu_a)_a\in\mathcal E$. At each index $a\in M/M_{1,j}$ with label $(v,c)$ and $j\sim i$, the allowed epimorphisms $L_v\to F_i$ are precisely the members of $[\mu_a]=\mu_a\Inn(F_i)$. Postcomposition gives a regular action of $\Inn(F_i)$ on $[\mu_a]$: it is transitive by definition, and if $\mu_a\alpha=\mu_a$ for $\alpha\in\Inn(F_i)$, then $\alpha=1$ because $\mu_a$ is surjective. Thus the choice of $\mu$ identifies $\mathcal E$ with the underlying set of
\begin{equation}\label{eq:epimorphism-family-product}
\prod_{i\in I}\prod_{j\sim i}\prod_{a\in M/M_{1,j}}(F_i/\mathbf{Z}(F_i)).
\end{equation}
For $y\in Y$, write $w$, $n_v$ and $q_c$ for its coordinates. By~\eqref{eq:general-array-action}, replacing $\mu_a$ by $\mu_a\Inn(f)$, where $f\in F_i$, replaces $(\mu^y)_{a^w}$ by $(\mu^y)_{a^w}\Inn(f^{\tau_i(q_c)})$. Hence the action of $Y$ on $\mathcal E$ is affine in the sense defined in Subsection~\ref{subsec:affine}: the translation part of $y$ is left multiplication by the element of~\eqref{eq:epimorphism-family-product} corresponding to $\mu^y$, and its linear part sends the factor indexed by $a$ to the factor indexed by $a^w$ via $f\mathbf{Z}(F_i)\mapsto f^{\tau_i(q_c)}\mathbf{Z}(F_i)$. If $y\in Y_0$, then its linear part is inner.

Apply Lemma~\ref{lem:refined-affine} with $(G,G_0)=(Y,Y_0)$, coefficient group $F_i/\mathbf{Z}(F_i)$ for each $j\sim i$, and a prime $p$ with $m_p=\lpp(m)$. Since $F_i\cong L/R_j$ for $j\sim i$, it follows from~\eqref{eq:general-incidence-energy} and Lemma~\ref{lem:prime-energy} that
\[
\begin{aligned}
\Psi&=\sum_{i\in I}\sum_{j\sim i}m_j\log(|F_i/\mathbf{Z}(F_i)|)\leq\sum_{j=1}^t m_j\log(|L/R_j|)=\Lambda\leq2(\ell-1)\log(2\ell),\\
\Psi_p&=\sum_{i\in I}\sum_{j\sim i}m_j\log(|F_i/\mathbf{Z}(F_i)|_p)\leq\sum_{j=1}^t m_j\log(|L/R_j|_p)=\Lambda_p\leq2(\ell-1).
\end{aligned}
\]
The number of $Y$-invariant $Y_0$-orbits on $\mathcal E$ is thus at most $2^{m\Psi/m_p+8m\Psi_p/\sqrt{\log(m_p)}}$, where
\[
\frac{m\Psi}{m_p}+\frac{8m\Psi_p}{\sqrt{\log(m_p)}}\leq\frac{2m(\ell-1)\log(2\ell)}{\lpp(m)}+\frac{16m(\ell-1)}{\sqrt{\log(\lpp(m))}}\leq\theta(m,\ell)
\]
by~\eqref{eq:theta}.

Finally, for $x\in S(\mu)$ and an index $a$ with label $(v,c)$, conjugation by $x_v$ induces conjugation by $\xi_c(x)$ under $\mu_a$. For the element $s\in S_0(\mu)$ corresponding to $x$, equation~\eqref{eq:general-array-action} therefore gives
\[
(\mu^s)_a=\Inn(x_v^{-1})\mu_a\Inn(\xi_c(x))=\mu_a.
\]
Hence $S_0(\mu)\leq(Y_0)_\mu$. For $y\in Y_\mu$, define $z=(z_v)_{v\in\overline{\Omega}}\in\prod_{v\in\overline{\Omega}}L_v$ by letting $z_{v^w}=n_v^{-1}x_v n_v$ for each $v\in\overline{\Omega}$. Since $\mu^y=\mu$, equation~\eqref{eq:general-array-action} gives
\[
\mu_{a^w}(z_{v^w})=\mu_a(x_v)^{\tau_i(q_c)}=\xi_c(x)^{\tau_i(q_c)}.
\]
Thus $z\in S(\mu)$ and $\xi_{c^w}(z)=\xi_c(x)^{\tau_i(q_c)}$. By the multiplication rule in $Y$, the conjugate $y^{-1}sy$ is the element of $S_0(\mu)$ corresponding to $z$. This implies that $S_0(\mu)\trianglelefteq Y_\mu$.
\end{proof}

\subsection{Counting complements}\label{subsec:abstract-complements}

For a fixed tuple and a family $\mu$ in a $Y$-invariant $Y_0$-orbit, we next count complements to $(Y_0)_\mu/S_0(\mu)$ in $Y_\mu/S_0(\mu)$. This count depends only on the abstract construction. Its application to orbital partitions will be made after the realization step in Subsection~\ref{subsec:realization}.

\begin{lemma}\label{lem:recovery-completions}
Fix a tuple in $\mathcal T_1$ or $\mathcal T_2$ and a representative $\mu\in\mathcal E$ of a $Y$-invariant $Y_0$-orbit, and suppose that $S(\mu)$ is subdirect in $\prod_{v\in\overline{\Omega}}L_v$. There exists an absolute constant $C$ such that the number of conjugacy classes of complements to $(Y_0)_\mu/S_0(\mu)$ in $Y_\mu/S_0(\mu)$ is at most
\begin{align}
2^{2\ell\log(|M|)+\frac{8m\ell}{\sqrt{\log(m)}}}&\quad\text{if the tuple lies in }\mathcal T_1,\label{eq:full-completions}\\
2^{Cd(M)m\sqrt\ell\,\log^{\frac{3}{2}}(\ell)}&\quad\text{if the tuple lies in }\mathcal T_2.\label{eq:outer-completions}
\end{align}
\end{lemma}

\begin{proof}
Write $S_0=S_0(\mu)$ and $A=(Y_0)_\mu/S_0$. Since the $Y_0$-orbit of $\mu$ is $Y$-invariant, the projection to the first coordinate maps $Y_\mu$ onto $M$. By Lemma~\ref{lem:compatible-quotient-maps}, there is an exact sequence
\[
1\longrightarrow A\longrightarrow Y_\mu/S_0\longrightarrow M\longrightarrow1.
\]
If this extension does not split, then there are no complements to count. We may therefore assume that it splits, in which case conjugacy and $A$-conjugacy of complements coincide.

We realize $A$ as an $M$-module for tuples in either $\mathcal T_1$ or $\mathcal T_2$. Let $Z_i=\mathbf{Z}(F_i)$, and let $V$ be the direct product with one copy of $Z_i$ at every $a\in M/M_{1,j}$ for $j\sim i$ and $i\in I$. Set $U=1$ for $\mathcal T_1$. For $\mathcal T_2$, let $U$ consist of the elements of $V$ whose entries agree at all indices with the same quotient label $c$.

For $y=(x,(q_c)_{i\in I,\ c\in M/Q_i})\in(Y_0)_\mu$, take $t_c=q_c\in F_i$ for $\mathcal T_1$, and for $\mathcal T_2$ take any $t_c\in F_i$ inducing $q_c\in\Inn(F_i)$. Since $y$ fixes $\mu$, equation~\eqref{eq:general-array-action} and the surjectivity of $\mu_a$ give $\Inn(t_c)=\Inn(\mu_a(x_v))$ at each index $a$ with label $(v,c)$. Hence
\begin{equation}\label{eq:general-central-difference}
t_c\mu_a(x_v)^{-1}\in Z_i.
\end{equation}
The assignment
\[
y\longmapsto\bigl(t_c\mu_a(x_v)^{-1}\bigr)_{i\in I,\ j\sim i,\ a\in M/M_{1,j}}U,
\]
where $a$ has label $(v,c)$, defines a homomorphism $(Y_0)_\mu\to V/U$. An element is in the kernel precisely when the values $\mu_a(x_v)$ agree at each $c$ and give the prescribed coordinate $q_c$ in $D_i$. Thus the kernel is $S_0$, and $A$ embeds in $V/U$.

If $y\in Y_\mu$ has first coordinate $g\in M$, then conjugation by $y$ takes the entry at $a$ to the entry at $a^g$ by $\tau_i(q_c(y))$, by~\eqref{eq:general-array-action}. This action on $V$ preserves $U$ and factors through $M$, since inner automorphisms act trivially on $Z_i$. The embedding of $A$ in $V/U$ is therefore $M$-equivariant. Once a complement is chosen, the $A$-conjugacy classes of complements are parametrized by $\mathrm H^1(M,A)$~\cite[Lemma~7]{Stix}. We bound its order using the above embedding.

For $\mathcal T_1$ we have $A\leq V$. At $a=M_{1,j}$ the action of $M_{1,j}$ on $Z_i$ is induced by $\beta_j$, by~\eqref{eq:compatible-outer-actions}. Hence
\[
V\cong\bigoplus_{i\in I}\bigoplus_{j\sim i}\Ind_{M_{1,j}}^M(Z_i).
\]
The cohomology sequence of $0\to A\to V\to V/A\to0$ gives $|\mathrm H^1(M,A)|\leq |(V/A)^M|\,|\mathrm H^1(M,V)|$. By Shapiro's lemma~\cite[Proposition~8]{Stix}, $\mathrm H^1(M,V)$ is the direct sum of the groups $\mathrm H^1(M_{1,j},Z_i)$. Since a $1$-cocycle is determined by its values on generators, Lemma~\ref{lem:prime-energy} gives
\[
|\mathrm H^1(M,V)|\leq\prod_{i\in I}\prod_{j\sim i}|Z_i|^{d(M_{1,j})}
\leq 2^{\log(|M|)\Lambda_0}\leq 2^{2\ell\log(|M|)}.
\]
Filtering the trivial section $(V/A)^M$ by the partial sums of the induced summands gives factors that are trivial sections of the individual summands. By Lemma~\ref{lem:trivial-sections},
\[
\log(|(V/A)^M|)\leq\sum_{i\in I}\sum_{j\sim i}\frac{2mm_j\log(|Z_i|)}{\sqrt{\pi\log(mm_j)}}
\leq\frac{2m\Lambda_0}{\sqrt{\pi\log(m)}}\leq\frac{8m\ell}{\sqrt{\log(m)}}.
\]
Multiplying these bounds gives~\eqref{eq:full-completions}.

For $\mathcal T_2$ we bound $|A|$ using the freedom to choose $t_c$ in~\eqref{eq:general-central-difference}. For each $c$, choose one index $b$ with label $(v_b,c)$ and take $t_c=\mu_b(x_{v_b})$. The entry at $b$ is then $1$, so this choice gives a unique representative of each image in $V/U$. Fix a real number $z>0$. At each $v$ take the image of $x_v$ under the product of the epimorphisms $\mu_a$ with block label $v$ and $\log(|F_i|)\leq z$. At the block labelled by $M_1$, the indices $a=M_{1,j}w$ with $w\in M_1$ give precisely the conjugates of $R_j$ under $N$, by~\eqref{eq:general-array-action}; each occurs once because $M_{1,j}$ is its stabilizer. An element of $Y$ with first coordinate $r_v$ transports these kernels to the block labelled by $v$. Thus, under the identification $L_v\cong L$, the kernels at each block run once through $\mathcal O(L,L_1)$, for every $\mu\in\mathcal E$. By Lemmas~\ref{lem:quotient-points} and~\ref{lem:joint-faithful}, this image has at most $2^{z\log(\ell)}$ possible values at each block. Over all $m$ blocks these values determine every normalized entry~\eqref{eq:general-central-difference} with $\log(|F_i|)\leq z$, including its reference value $\mu_b(x_{v_b})$.

There are at most $\Lambda/z$ remaining indices at each block, since the sum of $\log(|F_i|)$ over its indices is $\Lambda$. At each such index the normalized entry belongs to $Z_i$, of order at most $\ell$ by Lemma~\ref{lem:quotient-points}. Since the normalized entries determine the element of $A$, we obtain
\[
\log(|A|)\leq m(z+\Lambda/z)\log(\ell)\leq Cm\sqrt{\ell}\,\log^{\frac32}(\ell)
\]
on taking $z=\sqrt{\Lambda}$ and applying Lemma~\ref{lem:prime-energy}. A $1$-cocycle on $M$ is determined by its values on $d(M)$ generators, so $|\mathrm H^1(M,A)|\leq |A|^{d(M)}$. This gives~\eqref{eq:outer-completions}.
\end{proof}

\subsection{Realization by permutation groups}\label{subsec:realization}

We now return to the groups in Proposition~\ref{prop:general-recovery}. We prove that every such group gives data satisfying the abstract definitions, and that the resulting family recovers the equalities in Lemma~\ref{lem:block-pair-orbits}.

\begin{lemma}\label{lem:tuple-realization}
Let $\mathcal T\in\{\mathcal T_1,\mathcal T_2\}$. For every class in $\mathcal R$, there exist a tuple in $\mathcal T$, a representative $\mu$ of a $Y$-invariant $Y_0$-orbit on $\mathcal E$, and a subgroup $G\leq Y_\mu$ such that the following hold:
\begin{enumerate}[label=\textup{(\alph*)},leftmargin=*]
\item\label{item:compatible-realization} The projection to the first coordinate maps $G$ onto $M$, and $G$ gives rise to the given class in $\mathcal R$.
\item\label{item:compatible-subdirect} $S(\mu)$ is subdirect in $\prod_{v\in\overline{\Omega}}L_v$, and $G\cap Y_0\leq S_0(\mu)$.
\item\label{item:compatible-orbits} $G\cap Y_0$ and $S_0(\mu)$ have the same orbits on $(\overline{\Omega}\times\Delta)^2$.
\end{enumerate}
\end{lemma}

\begin{proof}
Fix a class in $\mathcal R$, and choose $G,K,\Omega$ and $\sigma$ giving rise to it. Identify the block $\Omega_{M_1}$ with $\Delta$ so that the group induced by $K$ is $L\in\mathcal L$, and let $\omega$ correspond to $1\in\Delta$. Set $H=G_\omega$, and let $K_0$ be the kernel of $K$ on this block. Its stabilizer in $G$ is $HK$, with $HK/K=M_1$. Since $HK$ normalizes $K$ and stabilizes the block, it normalizes $K_0$. The elements of $HK$ with a fixed image $x\in M_1$ under $HK\to HK/K=M_1$ induce a single coset of $L$ in $\mathbf{N}_{\Sym(\Delta)}(L)$. These cosets define a homomorphism $\rho\colon M_1\to\mathbf{N}_{\Sym(\Delta)}(L)/L$. The homomorphism $g\mapsto(gK,g|_\Delta)$ from $HK$ has kernel $K_0$ and image the group $N$ in~\eqref{eq:abstract-block-normalizer}, so we identify $HK/K_0$ with $N$ and $K/K_0$ with $L$.

Use the fixed orbit representatives $R_1,\ldots,R_t$ for this pair $(L,\rho)$, and let $C_j$ be the preimage of $R_j$ under $K\to L$. Write $i\sim j$ if $C_i$ and $C_j$ are $G$-conjugate, and let $I$ consist of the least index in each class. For $i\in I$, set
\[
Q_i=\mathbf{N}_G(C_i)/K,\quad E_i=\mathbf{N}_G(C_i)/C_i,\quad F_i=K/C_i=L/R_i.
\]
The subgroup $\mathbf{N}_{HK}(C_j)$ maps onto $M_{1,j}$ under $HK\to M_1$, so $Q_i$ contains $M_{1,i}$. The $G$-conjugates of $C_i$ are indexed by $M/Q_i$, with $g^{-1}C_i g$ labelled by $Q_i\overline g$.

For $j\sim i$, choose $g_j\in G$ with $g_j^{-1}C_ig_j=C_j$, taking $g_i=1$. By replacing $g_j$ within $\mathbf{N}_G(C_i)g_j$, arrange that $w_j=\overline{g_j}$ is the representative of $Q_i w_j$ in the fixed transversal. Define
\[
\varepsilon_j\colon L\longrightarrow F_i,\quad kK_0\longmapsto g_j k g_j^{-1}C_i.
\]
This is well-defined because $K_0\leq C_j$, and its kernel is $R_j$. Also $w_i=1$ and $\varepsilon_i$ is the natural quotient homomorphism modulo $R_i$. The inclusion $g_j\mathbf{N}_{HK}(C_j)g_j^{-1}\leq\mathbf{N}_G(C_i)$ gives $w_j M_{1,j} w_j^{-1}\leq Q_i$.

Let $\varphi_i$ be the outer action induced by the extension $E_i$ of $Q_i$ by $F_i$. To verify~\eqref{eq:compatible-outer-actions}, take $x\in M_{1,j}$ and choose $h\in\mathbf{N}_{HK}(C_j)$ with $hK=x$. The automorphism of $F_i$ induced by conjugation by $g_jhg_j^{-1}$ corresponds under $\varepsilon_j$ to the automorphism of $L/R_j$ induced by $h$. Their outer classes are $\varphi_i(w_j x w_j^{-1})$ and $\beta_j(x)$, respectively. The preimage of $M_{1,i}$ under $E_i\to Q_i$ is $\mathbf{N}_{HK}(C_i)/C_i$, which is identified with $\mathbf{N}_N(R_i)/R_i$ under $HK/K_0=N$. Thus all conditions defining $\mathcal T_1$ hold. Replacing each $E_i$ by its fixed equivalent representative preserves $F_i$, $Q_i$ and its outer action, and gives a tuple in $\mathcal T_1$. Keeping only the outer actions gives a tuple in $\mathcal T_2$.

Take the resulting tuple in $\mathcal T$, and construct $Y$ and $\mathcal E$ from it as above. Choose lifts $\gamma_v\in G$ of $r_v$, taking $\gamma_{M_1}=1$. The bijections $\delta\mapsto\delta^{\gamma_v}$ identify $\Omega$ with $\overline{\Omega}\times\Delta$ and its block actions with the copies $L_v$ used above. In these coordinates the action of $g\in G$ is
\[
(v,\delta)^g=(v^{\overline g},\delta^{n_v(g)}),\quad n_v(g)=\gamma_v g\gamma_{v^{\overline g}}^{-1}K_0\in N.
\]
For each $c\in M/Q_i$, choose a lift $g_c\in G$ of $d_c$. The element $g_c g g_{c^{\overline g}}^{-1}$ normalizes $C_i$. Its image modulo $C_i$, transported to the fixed representative $E_i$ if necessary, gives the coordinate $q_c(g)$ for $\mathcal T_1$; for $\mathcal T_2$, take $q_c(g)\in B_i$ to be the pair consisting of its image in $Q_i$ and its conjugation action on $K/C_i=F_i$. These coordinates respect multiplication in $Y$ and define an embedding of $G$ in $Y$, since its action on $\Omega$ is faithful. Under this embedding, $G\cap Y_0=K$.

The indices in~\eqref{eq:general-root-map} describe all occurrences of constituent kernels. Indeed, the stabilizer of $(\Delta,C_j)$ in $G$ has image $M_{1,j}$ in $M$ and contains $K$, so its translates are indexed by $M/M_{1,j}$. The translate indexed by $M_{1,j}w$ has block label $M_1w$ and kernel label $Q_i w_j w$. Thus the kernel at a label $c\in M/Q_i$ is $g_c^{-1}C_i g_c$. Define the quotient homomorphism
\[
\xi_c\colon K\longrightarrow F_i,\quad k\longmapsto g_c k g_c^{-1}C_i.
\]
At each index $a$ with label $(v,c)$, the kernel of $\xi_c$ contains $\Ker(\zeta_v)$, so there is a unique epimorphism $\mu_a\colon L_v\to F_i$ such that
\begin{equation}\label{eq:general-actual-array}
\mu_a(\zeta_v(k))=\xi_c(k)\quad\text{for every }k\in K.
\end{equation}
All constituent kernels on all blocks occur among these indices. Hence the equalities defining $S(\mu)$ are precisely those defining $S$ in Lemma~\ref{lem:block-pair-orbits}, and $S(\mu)=S$. The common-value homomorphisms $\xi_c$ on $S(\mu)$ extend those just defined on $K$.

At $a=M_{1,j}$, the lifts $g_c$ and $g_j$ of $w_j$ differ by an element of $K$, so the entry defined by~\eqref{eq:general-actual-array} lies in the class of $\varepsilon_j$. For $g\in G$ and $k\in K$, conjugation gives
\[
\zeta_{v^{\overline{g}}}(g^{-1}kg)=\zeta_v(k)^{n_v(g)},\quad
\xi_{c^{\overline{g}}}(g^{-1}kg)=\xi_c(k)^{\tau_i(q_c(g))}.
\]
By uniqueness in~\eqref{eq:general-actual-array}, \eqref{eq:general-array-action} sends $\mu_a$ to $\mu_{a^{\overline{g}}}$. Since $G$ maps onto $M$, every entry belongs to its prescribed class, so $\mu\in\mathcal E$ and $G\leq Y_\mu$. Thus $Y=Y_0G$ stabilizes the $Y_0$-orbit of $\mu$. Conjugating $G$ by an element of $Y_0$ takes $\mu$ to a chosen representative and preserves $\sigma$.

For this family, \eqref{eq:general-actual-array} gives $G\cap Y_0\leq S_0(\mu)$. Their images on $\Omega$ are $K$ and $S$, respectively, so they have the same ordered-pair orbits by Lemma~\ref{lem:block-pair-orbits}. Also $S$ is subdirect in $\prod_{v\in\overline{\Omega}}L_v$ since it contains $K$. These conclusions persist under the change of representative.
\end{proof}

\subsection{Proof of Proposition~\ref{prop:general-recovery}}

\begin{proof}[\rm\textbf{Proof of Proposition~\ref{prop:general-recovery}}]
Fix $\mathcal T\in\{\mathcal T_1,\mathcal T_2\}$. Lemma~\ref{lem:tuple-realization} gives a tuple in $\mathcal T$ and a $Y$-invariant $Y_0$-orbit on $\mathcal E$ for every class in $\mathcal R$. Choose one representative $\mu$ from each such orbit. If Lemma~\ref{lem:tuple-realization} gives $G\leq Y_\nu$ with $\nu$ in the orbit represented by $\mu$, then choose $y_0\in Y_0$ with $\nu^{y_0}=\mu$, so that $y_0^{-1}Gy_0\leq Y_\mu$. The element $y_0$ preserves the block labels given by $\sigma$, so this conjugate represents the same class in $\mathcal R$.

For a representative that occurs, write $S_0=S_0(\mu)$ and $A=(Y_0)_\mu/S_0$. The subgroup $G\leq Y_\mu$ supplied by Lemma~\ref{lem:tuple-realization} maps onto $M$ and satisfies $G\cap Y_0\leq S_0$, and the groups $G\cap Y_0$ and $S_0$ have the same orbits on $(\overline{\Omega}\times\Delta)^2$. Hence $GS_0/S_0$ is a complement to $A$ in $Y_\mu/S_0$. Given an ordered pair and $s\in S_0$, an element of $G\cap Y_0$ has the same effect as $s$ on that pair, so $GS_0$ and $G$ have the same orbital partition. The full preimage in $Y_\mu$ of the complement is $GS_0$, and $A$-conjugate complements have preimages conjugate by an element of $(Y_0)_\mu$, whose action preserves $\sigma$. Thus each $A$-conjugacy class of complements accounts for at most one class in $\mathcal R$ for this tuple and this representative.

For $\mathcal T_1$, Lemma~\ref{lem:recovery-tuples} bounds the number of tuples by~\eqref{eq:general-full-front}, Lemma~\ref{lem:compatible-quotient-maps} bounds the number of invariant orbits for each tuple by $2^{\theta(m,\ell)}$, and Lemma~\ref{lem:recovery-completions} bounds the number of complement classes for each representative by~\eqref{eq:full-completions}. Multiplying these bounds and taking logarithms gives~\eqref{eq:general-full}, after absorbing $2\ell\log(|M|)$ into its first term and enlarging the absolute constant $C$.

For $\mathcal T_2$, multiply~\eqref{eq:general-outer-front}, the bound $2^{\theta(m,\ell)}$ for the invariant orbits, and~\eqref{eq:outer-completions}. Their product gives~\eqref{eq:general-outer}.
\end{proof}

\section{Minimal normal subgroups with prescribed orbit size}\label{sec:minimal-normal}

In this section, we bound the orders and the number of permutational isomorphism classes of minimally transitive permutation groups having a minimal normal subgroup whose orbits have prescribed size.

\begin{lemma}\label{lem:simple-power}
Let $T$ be a nonabelian finite simple group, let $r\geq1$ be an integer, and let $H\leq T^r$ be a core-free subgroup.
Set $b=[T^r:H]$. Then
\[
|T^r|\leq b^{b^{\frac{1}{r}}}\quad\text{and}\quad r\leq\frac{\log(b)}{\log(5)}.
\]
\end{lemma}

\begin{proof}
Let $a$ be the minimum index of a proper subgroup of $T$. Then $|T|\leq a!$ and $a\geq5$.
By~\cite[Proposition~5.2.7(ii) and the observation following its proof]{KL}, we have $b\geq a^r$, and so
\[
|T^r|\leq(a!)^r\leq(a^r)^a\leq b^{b^{\frac{1}{r}}}.
\]
The inequality $r\leq\log(b)/\log(5)$ follows from $b\geq a^r\geq5^r$.
\end{proof}

\begin{lemma}\label{lem:minimal-normal-local}
For each integer $b\geq2$, let $\mathcal L(b)$ be the set of $\Sym(b)$-conjugacy classes of subgroups $L\leq\Sym(b)$ that are permutationally isomorphic to the group induced on an orbit of size $b$ by a minimal normal subgroup of a minimally transitive permutation group.
Then
\[
|\mathcal L(b)|\leq2^{O(\sqrt b\log^2(b))}\quad\text{and}\quad|\mathbf{N}_{\Sym(b)}(L)|=2^{O(\sqrt b\log(b))}
\]
for every $L$ represented in $\mathcal L(b)$.
\end{lemma}

\begin{proof}
Fix $L$ represented in $\mathcal L(b)$, let $N=\mathbf{N}_{\Sym(b)}(L)$, and choose a minimally transitive permutation group $G\leq\Sym(\Omega)$ and a minimal normal subgroup $X$ of $G$ inducing $L$ on an orbit in $\Omega$ of size $b$.
If $X$ is abelian, then $L=C_p^a$ is regular, where $p$ is a prime and $a\geq1$ is an integer such that $b=p^a$. In this case, $N=\AGL_a(p)$  with $|N|\leq p^{a^2+a}=2^{O(\log^2(b))}$, and all such $L$ lie in one $\Sym(b)$-conjugacy class.

Suppose that $X=T^t$ for a nonabelian simple group $T$.
Then $L=T^r$ for some $r\in\{1,\ldots,t\}$.
If $r\geq2$, then Lemma~\ref{lem:simple-power} gives $|L|\leq b^{b^{1/r}}\leq b^{\sqrt b}$.
Assume that $r=1$. Let $H<T$ be a point stabilizer in $T$, and choose a maximal subgroup $M$ of $T$ containing $H$.
Let $a=[T:M]\leq b$.
The action on $T/M$ is faithful and primitive.
If this primitive group does not contain $\Alt(a)$, then Mar\'oti's bound~\cite[Corollary~1.1(ii)]{Maroti} gives $|T|<50a^{\sqrt a}\leq2^{O(\sqrt b\log(b))}$. Otherwise, $T=\Alt(a)$ and $M=\Alt(a-1)$ in that action.
Increasing the implied absolute constant if necessary, we may assume that $a\geq7$.
If $H<M$, then $[M:H]\geq a-1$, and so $b\geq a(a-1)$, which still yields $|T|=2^{O(a\log(a))}=2^{O(\sqrt b\log(b))}$.

We exclude the possibility that $T=\Alt(a)$ and $H=M=\Alt(a-1)$. Suppose for a contradiction that this is the case.
For prime $b$, take $V=\AGL_1(b)\cap\Alt(b)$.
For $b=ef$ with integers $e\geq2$ and $f\geq2$, take $V=(\Sym(e)\wr\Sym(f))\cap\Alt(b)$ in the imprimitive action.
In both cases, $V$ is a proper transitive subgroup of $\Alt(b)$ whose $\Alt(b)$-conjugacy class is invariant under $\Aut(\Alt(b))=\Sym(b)$.
Now let $W=V^t\leq X=\Alt(b)^t$.
Since $r=1$, on each $X$-orbit in $\Omega$, exactly one simple direct factor of $X$ acts nontrivially, and its induced action is the natural action of $\Alt(b)$ up to an automorphism of $\Alt(b)$.
Thus, $W$ is transitive on every $X$-orbit.
Note that conjugation by $G$ permutes the simple direct factors of $X$ and acts on each factor through an automorphism.
We conclude that every $G$-conjugate of $W$ is $X$-conjugate to it. Hence
\[
G=X\mathbf{N}_G(W),
\]
and so $\mathbf{N}_G(W)$ is transitive on the set of $X$-orbits, as $G$ is. Since $\mathbf{N}_G(W)$ contains $W$, which is transitive on every $X$-orbit, it follows that $\mathbf{N}_G(W)$ is transitive on $\Omega$, contradicting minimal transitivity of $G$.

So far we have proved $|L|=2^{O(\sqrt b\log(b))}$. The centralizer of the transitive subgroup $L$ in $\Sym(b)$ has order at most $b$, and $\Aut(L)=\Aut(T^r)=\Aut(T)\wr\Sym(r)$. Then since $|\Out(T)|<|T|^{1/3}$ (see for example~\cite[Lemma~4.8]{Fawcett}) and $r\leq\log(b)/\log(5)$ by Lemma~\ref{lem:simple-power},
\[
|N|\leq b|\Aut(T^r)|\leq b|T^r|^{\frac{4}{3}}r!=b|L|^{\frac{4}{3}}r!=2^{O(\sqrt b\log(b))},
\]
as required.

To bound $|\mathcal L(b)|$, let $c=C\sqrt b\log(b)$ for a sufficiently large absolute constant $C$, so that the possible $T^r$ have order at most $2^c$.
By~\cite[Theorem~6.1]{KLST}, there are at most two isomorphism classes of nonabelian characteristically simple groups of any given order, and hence at most $2^{c+1}$ isomorphism classes of possible $T^r$.
For each abstract $T^r$, there are at most
\[
|T^r|^{2\log(b)}\leq2^{2c\log(b)}
\]
subgroups of index $b$ by~\cite[Corollary~1.7.5]{LS}. It follows that
\[
|\mathcal L(b)|\leq2^{c+1}\cdot2^{2c\log(b)}=2^{O(c\log(b))}=2^{O(\sqrt b\log^2(b))},
\]
completing the proof.
\end{proof}

We are now ready to establish the main result of this section.

\begin{proposition}\label{prop:large-orbit-catalogue}
Let $2\leq b\leq m$ be integers such that $b$ divides $m$, and let $\mathcal M(m,b)$ be the set of $\Sym(m)$-conjugacy classes of minimally transitive subgroups $M\leq\Sym(m)$ that have a minimal normal subgroup with orbit size $b$. Then
\[
|\mathcal M(m,b)|\leq2^{O\left(\frac{m\log^2(m)}{\sqrt b}+\sqrt b\log^2(b)\right)}\quad\text{and}\quad|M|=2^{O\left(\frac{m\log(m)}{\sqrt b}\right)}
\]
for every group $M$ represented in $\mathcal M(m,b)$.
\end{proposition}

\begin{proof}
Fix a representative $M$ of a class in $\mathcal M(m,b)$ and a minimal normal subgroup $X$ of $M$ with orbit size $b$.
This minimal normal subgroup has $k=m/b$ orbits. Let $\Delta$ be one of these orbits, and let $L\leq\Sym(b)$ be the group induced by $X$ on $\Delta$, after identifying $\Delta$ with $\{1,\ldots,b\}$.
Set $N=\mathbf{N}_{\Sym(b)}(L)$, so that $L$ represents a class in $\mathcal L(b)$ from Lemma~\ref{lem:minimal-normal-local}.
Since $X\trianglelefteq M$, the group induced by the stabilizer $M_\Delta$ on $\Delta$ normalizes $L$ and hence is contained in $N$.
Thus we may assume $M\leq N\wr\Sym(k)$ in its imprimitive action on $m$ points. By Lemma~\ref{lem:minimal-normal-local},
\[
|N\wr\Sym(k)|=2^{O(k\sqrt b\log(b)+k\log(k))}=2^{O\left(\frac{m\log(m)}{\sqrt b}\right)},
\]
and hence $|M|=2^{O(m\log(m)/\sqrt{b})}$. For each fixed representative $L$, Lemma~\ref{lem:rank-generators}\ref{item:generators} shows that each minimally transitive subgroup of $N\wr\Sym(k)$ is generated by at most $\log(m)$ elements, and so the number of such subgroups is at most
\[
|N\wr\Sym(k)|^{\log(m)}=2^{O\left(\frac{m\log^2(m)}{\sqrt b}\right)}.
\]
Then multiplying by the bound $2^{O(\sqrt b\log^2(b))}$ for $|\mathcal L(b)|$ from Lemma~\ref{lem:minimal-normal-local} gives the bound for $|\mathcal M(m,b)|$.
\end{proof}

\section{Proof of Theorem~\ref{thm:main}}\label{sec:counting}

Let $\mathcal{P}(n)$ be the set of isomorphism classes of orbital partitions of minimally transitive permutation groups of degree $n$.
The key step in the proof of Theorem~\ref{thm:main} is the following proposition, which gives a subexponential bound for $|\mathcal{P}(n)|$.

\begin{proposition}\label{prop:partitions}
There exists an absolute constant $C$ such that, for every integer $n\geq3$,
\begin{equation}\label{eq:partition-count}
|\mathcal{P}(n)|\leq2^{\frac{Cn}{\sqrt{\log(\log(n))}}}.
\end{equation}
\end{proposition}

\begin{proof}
After increasing $C$ to cover the finitely many smaller values of $n$, it suffices to consider sufficiently large $n$, so that $\log^{24}(n)<n$.
Let $G$ be a minimally transitive permutation group of degree $n$.
Starting with $K_0=1$, construct normal subgroups $K_i$ of $G$ recursively, writing $k_i$ for the common size of the $K_i$-orbits.
If $k_i<n$, then choose a minimal normal subgroup $N_i$ of $G/K_i$ acting on the set of $K_i$-orbits, and let $K_{i+1}$ be the preimage under the natural homomorphism $G\to G/K_i$ of the kernel of the action of $G/K_i$ on the set of $N_i$-orbits.
Each $K_i$ is the kernel of the action on its own orbit set, the induced group $G/K_i$ is minimally transitive, and the common size of the $N_i$-orbits is $k_{i+1}/k_i\geq2$.
The orbit sizes therefore strictly increase from $k_0=1$ to $n$. Let $j$ be the least index such that $k_j\geq\log^6(n)$.
Then $k_{j-1}<\log^6(n)$.
We bound the number of isomorphism classes of orbital partitions in each of the following two cases by the right-hand side of~\eqref{eq:partition-count}. Since every element of $\mathcal{P}(n)$ arises in at least one of the two cases, adding these two bounds proves~\eqref{eq:partition-count} after increasing $C$.

\smallskip\noindent
\textsf{Case~1: $k_j\leq\log^{24}(n)$.}
Let $\ell=k_j$, $m=n/\ell$, $K=K_j$ and $M=G/K$. Then
\begin{equation}\label{eq:orbit-size-interval}
\log^6(n)\leq \ell\leq\log^{24}(n).
\end{equation}
The quotient $M$ is minimally transitive on the set of $m$ orbits of $K$, and we identify this set with $\{1,\ldots,m\}$.
By Lemma~\ref{lem:rank-generators}\ref{item:generators}, $d(M)\leq\log(m)\leq\log(n)$.
Counting the possible generating tuples, we see that the number of possibilities for $M\leq\Sym(m)$ is at most
\[
(m!)^{\log(n)}\leq2^{m\log^2(n)}=2^{\frac{n\log^2(n)}{\ell}}\leq2^{\frac{n}{\log^4(n)}}.
\]
For each such $M$, let $M_1$ be the stabilizer of $1$ in $M$. Then $d(M_1)\leq m$ by Lemma~\ref{lem:general-generators}.

We apply~\eqref{eq:general-outer} with quotient group $M$, point stabilizer $M_1$, normal subgroup $K$, block size $\ell$, and $\mathcal L$ taken to be a set of representatives of the permutational isomorphism classes of transitive permutation groups of degree $\ell$.
By~\eqref{eq:orbit-size-interval}, $\ell^2\log^2(\ell)\log(n)+\ell\log^2(n)=O(\log^{49}(n)\log^2(\log(n)))=o(n/\sqrt{\log(\log(n))})$.
Moreover, Lemma~\ref{lem:rank-generators}\ref{item:normalizer} together with $d(M_1)\leq m$ gives
\[
d(M_1)\max_{L\in\mathcal L}\log([\mathbf{N}_{\Sym(\ell)}(L):L])
\leq\frac{2n}{\sqrt{\log(\ell)}}+\frac{n\log^2(\ell)}{\ell}=O\left(\frac{n}{\sqrt{\log(\log(n))}}\right),
\]
and by Lemma~\ref{lem:uniform-theta}, $\theta(m,\ell)=O(n/\sqrt{\log(\log(n))})$. In view of $d(M)\leq\log(n)$, $m=n/\ell$ and the bounds in~\eqref{eq:orbit-size-interval}, we also have
\[
\begin{aligned}
d(M)m\sqrt\ell\,\log^{\frac{3}{2}}(\ell)&\leq\frac{n\log(n)\log^{\frac{3}{2}}(\ell)}{\sqrt\ell}\\
&\leq\frac{n\bigl(24\log(\log(n))\bigr)^{\frac{3}{2}}}{\log^2(n)}=O\left(\frac{n}{\sqrt{\log(\log(n))}}\right).
\end{aligned}
\]
Consequently, for each fixed $\ell$ and $M$, \eqref{eq:general-outer} gives at most $2^{O(n/\sqrt{\log(\log(n))})}$ isomorphism classes of orbital partitions arising from these groups $G$.
Multiplying by $2^{n/\log^4(n)}$, the bound for the number of possibilities for $M$, and summing over at most $n$ values of $\ell$ gives at most $2^{O(n/\sqrt{\log(\log(n))})}$ isomorphism classes of orbital partitions in Case~1.

\smallskip\noindent
\textsf{Case~2: $k_j>\log^{24}(n)$.}
Let $\ell=k_{j-1}<\log^6(n)$, $m=n/\ell$, $K=K_{j-1}$ and $M=G/K$.
The quotient $M$ is a minimally transitive permutation group of degree $m$.
It has a minimal normal subgroup with orbit size
\[
b=k_j/\ell>\log^{24}(n)/\ell>\log^{18}(n).
\]
For fixed $\ell$ and $b$, Proposition~\ref{prop:large-orbit-catalogue} gives at most $2^{O(n/\sqrt{\log(\log(n))})}$ possibilities for $M$, up to permutational isomorphism, since $m=n/\ell\leq n$ and
\[
\frac{m\log^2(m)}{\sqrt b}+\sqrt b\log^2(b)\leq\frac{n\log^2(n)}{\ell\sqrt b}+\sqrt b\log^2(b)\leq\frac{n}{\log^7(n)}+\sqrt n\,\log^2(n)
=O\left(\frac{n}{\sqrt{\log(\log(n))}}\right).
\]
It also follows from Proposition~\ref{prop:large-orbit-catalogue} that
\begin{equation}\label{eq:jump-order}
\log(|M|)\leq \frac{C_0n\log(n)}{\ell\sqrt b}
\end{equation}
for some absolute constant $C_0$.
If $\ell=1$, then $K=1$ and $G=M$, in which case Proposition~\ref{prop:large-orbit-catalogue} gives at most $2^{O(n/\sqrt{\log(\log(n))})}$ isomorphism classes of orbital partitions, as permutationally isomorphic groups have isomorphic orbital partitions.

Assume $\ell\geq2$ and apply~\eqref{eq:general-full} with quotient group $M$, a point stabilizer of $M$, normal subgroup $K$, and $\mathcal L$ taken to be a set of representatives of the permutational isomorphism classes of transitive permutation groups of degree $\ell$.
Since $\ell<\log^6(n)$, we derive from~\eqref{eq:jump-order} that
\[
\begin{aligned}
\ell^2\log^2(\ell)\log(|M|)&\leq \frac{C_0n\ell\log(n)\log^2(\ell)}{\sqrt b}=O(n/\sqrt{\log(\log(n))}),\\
\ell\log(\ell)(\log(n)+1)\log(|M|)&\leq \frac{C_0n\log(\ell)(\log(n)+1)\log(n)}{\sqrt b}=O(n/\sqrt{\log(\log(n))}).
\end{aligned}
\]
Also $\ell\log^2(n)<\log^8(n)$, and $\lpp(m)\leq m$ gives $8m\ell/\sqrt{\log(m)}\leq\theta(m,\ell)/2$.
Then~\eqref{eq:general-full} together with Lemma~\ref{lem:uniform-theta} gives at most $2^{O(n/\sqrt{\log(\log(n))})}$ isomorphism classes of orbital partitions for each fixed $\ell$ and quotient action $M$. Since there are $2^{O(n/\sqrt{\log(\log(n))})}$ possibilities for $M$, up to permutational isomorphism, for each fixed pair $(\ell,b)$, while there are at most $n^2$ pairs $(\ell,b)$, the number of isomorphism classes of orbital partitions in this case is thus at most
\[
n^2\cdot2^{O\left(\frac{n}{\sqrt{\log(\log(n))}}\right)}\cdot2^{O\left(\frac{n}{\sqrt{\log(\log(n))}}\right)}
=2^{O\left(\frac{n}{\sqrt{\log(\log(n))}}\right)}.\qedhere
\]
\end{proof}

For a transitive permutation group $G$ on a set $\Omega$, let $q(G)$ be the number of orbits of $G$ on the set of $2$-element subsets of $\Omega$. Then $q(G)$ depends only on the orbital partition of $G$. The numbers of invariant digraphs and graphs on $\Omega$ are $2^{r(G)-1}$ and $2^{q(G)}$, respectively.
For each integer $n\geq2$, define $\beta(n)$, $\gamma(n)$ and $R_n$ by the following table, where $C_n$ denotes the cyclic group of order $n$, and $D_n$ the dihedral group of order $n$.
\begin{equation}\label{eq:counting-data}
\begin{array}{c|c|c|c}
n&\beta(n)&\gamma(n)&R_n\\ \hline
\text{odd}&11/24&1/2&C_n\\
\text{even but not a $2$-power}&11/16&3/4&D_n\\
2^f&11/16&1&C_2^f
\end{array}
\end{equation}

\begin{lemma}\label{lem:graph-rank}
Let $G$ be a nonregular minimally transitive permutation group of degree $n\geq2$.
Then $q(G)\leq\beta(n)n$.
\end{lemma}

\begin{proof}
The proof of~\cite[Corollary~2]{ETXZ} shows that $G$ preserves at most $2^{11n/16}$ graphs.
Since the number of these graphs is $2^{q(G)}$, the result follows when $n$ is even.
For the rest of the proof, assume now that $n$ is odd.
Let $H$ be a point stabilizer in $G$.

Suppose that all subdegrees of $G$ are at most two.
By Bergman and Lenstra~\cite[Theorem~1]{BL}, in the form stated in~\cite[Section~1]{Spiga}, one of the following holds, as $G$ is nonregular:
\begin{enumerate}[label=\textup{(\roman*)},leftmargin=*]
\item\label{item:small-subdegrees-order-two} $|H|=2$;
\item\label{item:small-subdegrees-normal-two} $G$ has a normal elementary abelian $2$-subgroup $N$ containing $H$ with $[N:H]=2$.
\end{enumerate}
In case~\ref{item:small-subdegrees-order-two}, the nonidentity element of $H$ has sign $(-1)^n=-1$ in the regular action of $G$ on itself.
The kernel of this sign homomorphism therefore has order $n$ and intersects $H$ trivially, giving a proper transitive subgroup of $G$, contrary to its minimal transitivity.
In case~\ref{item:small-subdegrees-normal-two}, the $N$-orbits have size two, contrary to $n$ being odd.

Consequently, there is a suborbit of size at least three.
Spiga~\cite[Theorem~1.1]{Spiga} shows that the union $B$ of these suborbits has size at least $n/6$.
Each $G$-orbit on the set of $2$-element subsets of the point set of $G$ is the edge set of a graph.
At the point stabilized by $H$, the neighbour set in each such graph is either a single suborbit or the union of two paired suborbits.
Since paired suborbits have equal size, each of these neighbour sets is either contained in $B$ or disjoint from $B$.
Moreover, each graph has even valency by the handshaking lemma, since $n$ is odd. Hence its valency is at least two, and at least four if its neighbour set is contained in $B$. This implies that
\[
q(G)\leq\frac{n-1-|B|}{2}+\frac{|B|}{4}\leq\frac{11n}{24}-\frac12.\qedhere
\]
\end{proof}

\begin{lemma}\label{lem:cayley-lower-bounds}
For all sufficiently large $n$,
\[
|\mathcal C(n)|\geq2^{n-\log^2(n)-2}\quad\text{and}\quad
|\mathcal C(n)\cap\mathcal U(n)|\geq2^{\gamma(n)n-\log^2(n)-2}.
\]
\end{lemma}

\begin{proof}
The proof of Morris and Spiga~\cite[Theorem~1.5]{MS} gives at least $2^{n-\log^2(n)-1}$ isomorphism classes of digraphical regular representations of $C_n$, for sufficiently large $n$, when loops are allowed.
%Adding or removing all loops preserves the automorphism group and pairs the isomorphism classes with and without loops.
Thus excluding loops contributes a factor of $1/2$, giving the bound for $\mathcal C(n)$.

For graphs, take $R_n$ from~\eqref{eq:counting-data}.
The proofs cited below give at least $2^{(n+i(R_n))/2-\log^2(n)-2}$ isomorphism classes of graphs on $R_n$ with full automorphism group of order $n$ or $2n$, where the factor of $1/2$ for excluding loops is included in this bound.
If $n$ is odd or a $2$-power, then $R_n$ is abelian, and we use the proof of Dobson, Spiga and Verret~\cite[Theorem~6.2]{DSV}.
The graphs counted there have full automorphism group $R_n\rtimes\langle\iota\rangle$, where $\iota$ is inversion, of order $2n$ for odd $n$ and $n$ for a $2$-power.
If $n$ is even but not a $2$-power, then $R_n=D_n$, and we use the proof of Xia and Zheng~\cite[Corollary~1.4]{XZGRR}, which counts graphical regular representations.
Its hypotheses hold because $D_n$ is nonabelian and $i(D_n)=n/2+i(C_{n/2})>n/2$, whereas in a generalized dicyclic group every element outside its distinguished abelian subgroup of index two has order four.
Since $(n+i(R_n))/2\geq\gamma(n)n$ in all cases, the bound for $|\mathcal C(n)\cap\mathcal U(n)|$ follows.
\end{proof}

\begin{proof}[\rm\textbf{Proof of Theorem~\ref{thm:main}}]
Every non-Cayley vertex-transitive digraph or graph admits a nonregular minimally transitive subgroup of automorphisms.
For each class in $\mathcal{P}(n)$ that contains the orbital partition of a nonregular minimally transitive group, choose one such group $G$.
The arc sets of the digraphs admitting $G$ as a group of automorphisms are precisely the unions of its nondiagonal orbitals; for graphs, these unions must also be symmetric.
An isomorphism of orbital partitions carries each such union to the arc set of an isomorphic digraph or graph, so every non-Cayley vertex-transitive digraph or graph is isomorphic to one admitting a chosen group $G$.
For each chosen $G$, there are $2^{r(G)-1}$ such digraphs and $2^{q(G)}$ such graphs.
By Proposition~\ref{prop:partitions}, Lemma~\ref{lem:rank-generators}\ref{item:rank} and Lemma~\ref{lem:graph-rank},
% \begin{equation}\label{eq:non-cayley-counts}
\[
|\mathcal V(n)\setminus\mathcal C(n)|\leq2^{\frac{3n}{4}+F(n)}\quad\text{and}\quad|\mathcal U(n)\setminus\mathcal C(n)|\leq2^{\beta(n)n+F(n)},
\]
% \end{equation}
where $F(n)=Cn/\sqrt{\log(\log(n))}=o(n)$.
Let $E(n)=F(n)+\log^2(n)+2=o(n)$.
Together with Lemma~\ref{lem:cayley-lower-bounds}, these bounds give
\[
\begin{aligned}
\frac{|\mathcal C(n)|}{|\mathcal V(n)|}&=1-\frac{|\mathcal V(n)\setminus\mathcal C(n)|}{|\mathcal V(n)|}\geq1-\frac{|\mathcal V(n)\setminus\mathcal C(n)|}{|\mathcal C(n)|}\geq1-2^{-\frac{n}{4}+E(n)},\\
\frac{|\mathcal C(n)\cap\mathcal U(n)|}{|\mathcal U(n)|}&=1-\frac{|\mathcal U(n)\setminus\mathcal C(n)|}{|\mathcal U(n)|}\geq1-\frac{|\mathcal U(n)\setminus\mathcal C(n)|}{|\mathcal C(n)\cap\mathcal U(n)|}\geq1-2^{-(\gamma(n)-\beta(n))n+E(n)}.
\end{aligned}
\]
By~\eqref{eq:counting-data}, we have $\gamma(n)-\beta(n)\geq1/24$. Hence the theorem follows.
\end{proof}

\begin{remark}\label{rem:even-proportions}
For even $n$, the values of $\gamma(n)-\beta(n)$ are $1/16$ when $n$ is not a $2$-power and $5/16$ when $n$ is a $2$-power.
The proof of Theorem~\ref{thm:main} therefore gives the stronger bounds
\[
\frac{|\mathcal C(n)\cap\mathcal U(n)|}{|\mathcal U(n)|}\geq
\begin{cases}
1-2^{-\frac{n}{16}+o(n)}&\text{if $n$ is even but not a $2$-power}\\
1-2^{-\frac{5n}{16}+o(n)}&\text{if $n$ is a $2$-power}
\end{cases}
\]
as the even integer $n\to\infty$ through the respective subsequences.
\end{remark}

\section*{AI disclosure}

This work arose from collaboration between the authors and GPT-6.1 (OpenAI), which was used to explore possible arguments and produce intermediate drafts. The authors organized the investigation in stages: first considering prime-power degree, then solvable groups, and finally arbitrary groups. This approach divided the problem into successive tasks, with intermediate results developed before proceeding to the next stage.

The proof strategy and arguments in Sections~\ref{sec:minimal-normal} and~\ref{sec:counting}, using minimal normal subgroups, their orbit partitions and wreath product embeddings, are from the second author's earlier unpublished work on counting minimally transitive permutation groups. In the present paper, these arguments are combined with Proposition~\ref{prop:general-recovery} to complete the proof.

The principal mathematical contribution of AI was the discovery of Proposition~\ref{prop:general-recovery} and the supporting arguments developed in Sections~\ref{sec:auxiliary} and~\ref{sec:recovery}.
The authors have critically examined the AI-produced arguments, rewritten them, and substantially simplified a few proofs.

\end{document}